\documentclass[12pt]{amsart}
\usepackage{amsmath,amssymb,amscd,amsthm,latexsym}
\usepackage{euscript}
\DeclareMathAlphabet{\EuRm}{U}{eur}{m}{n}
\SetMathAlphabet{\EuRm}{bold}{U}{eur}{b}{n}
\begin{document} 
%
%
\swapnumbers
\newtheorem{thm}{Theorem}[section]
\newtheorem{lemma}[thm]{Lemma}
\newtheorem{prop}[thm]{Proposition}
\newtheorem{fact}[thm]{Fact}
\newtheorem{cor}[thm]{Corollary}
\theoremstyle{definition}
\newtheorem{defn}[thm]{Definition}
\newtheorem{example}[thm]{Example}
\newtheorem{examples}[thm]{Examples}
\newtheorem{claim}[thm]{Claim}
\newtheorem{con}[thm]{Convention}
\newtheorem{obs}[thm]{Observation}

\theoremstyle{remark}
\newtheorem{assume}[thm]{Assumption}
\newtheorem{remark}[thm]{Remark}
\newtheorem{note}[thm]{Note}
\newtheorem{notation}[thm]{Notation}
\newtheorem{aside}[thm]{Aside}
\newtheorem{ack}[thm]{Acknowledgements}
\numberwithin{equation}{section}
\numberwithin{figure}{section}
%
%
\newcommand{\sect}{\setcounter{thm}{0}\section}
%
%
\newcommand{\xra}[1]{\xrightarrow{#1}}
\newcommand{\xla}[1]{\xleftarrow{#1}}
\newcommand{\xsim}{\xrightarrow{\sim}}
\newcommand{\hra}{\hookrightarrow}
\newcommand{\epic}{\to\hspace{-5 mm}\to}
\newcommand{\adj}[2]{\substack{{#1}\\ \rightleftharpoons \\ {#2}}}
\newcommand{\ccsub}[1]{\circ_{#1}}
\newcommand{\DEF}{:=}
\newcommand{\EQUIV}{\Leftrightarrow}
\newcommand{\hsp}{\hspace{10 mm}}
\newcommand{\hs}{\hspace{5 mm}}
\newcommand{\hsm}{\hspace{3 mm}}
\newcommand{\vsm}{\vspace{2 mm}}
\newcommand{\vs}{\vspace{5 mm}}
\newcommand{\rest}[1]{\lvert_{#1}}
\newcommand{\lra}[1]{\langle{#1}\rangle}
\newcommand{\llrr}[1]{\langle\!\langle{#1}\rangle\!\rangle}
\newcommand{\llrra}[1]{\overline{\langle\!\langle{#1}\rangle\!\rangle}}
\newcommand{\lamt}[3]{\lambda_{3}({#1}\otimes{#2}\otimes{#3})}
\newcommand{\lamf}[4]{\lambda_{4}({#1}\otimes{#2}\otimes{#3}\otimes{#4})}
%
%
\newcommand{\pis}{\pi_{\ast}}
\newcommand{\pin}[2]{\hat{\pi}_{{#1}}{#2}}
%
%
\newcommand{\CC}{\mathbb C}
\newcommand{\bL}[1]{{\mathbb L}\lra{#1}}
\newcommand{\NN}{\mathbb N}
\newcommand{\QQ}{\mathbb Q}
\newcommand{\ZZ}{\mathbb Z}
\newcommand{\bn}{\mathbf{n}}
%
%
\newcommand{\Aa}{{\mathcal A}}
\newcommand{\Ab}{{\EuScript Ab}}
\newcommand{\Abgp}{{\EuScript AbGp}}
\newcommand{\Alg}[1]{{#1}\text{-}{\EuScript Alg}}
\newcommand{\C}{{\mathcal C}}
\newcommand{\Cat}{\textit{$\EuScript C$\hskip -.85pt at}}
\newcommand{\D}{{\mathcal D}}
\newcommand{\DGA}{{\EuScript DGA}}
\newcommand{\DGAo}{\DGA_{1}}
\newcommand{\DGC}{{\EuScript DGC}}
\newcommand{\DGL}{{\EuScript DGL}}
\newcommand{\DGLz}{\DGL_{0}}
\newcommand{\DGN}{{\EuScript DGN}}
\newcommand{\E}{{\mathcal E}}
\newcommand{\F}{{\mathcal F}}
\newcommand{\G}{{\mathcal G}}
\newcommand{\Lie}{{\EuScript Lie}}
\newcommand{\LL}{{\mathcal L}}
\newcommand{\M}{{\mathcal M}}
\newcommand{\Ner}{{\mathcal N}}
\newcommand{\OO}{{\mathcal O}}
\newcommand{\PP}{{\mathcal P}}
\newcommand{\Pa}{$\Pi$-algebra}
\newcommand{\PAlg}{\Alg{\Pi}}
\newcommand{\R}{{\mathcal R}}
\newcommand{\Ss}{{\mathcal S}}
\newcommand{\Sa}{\Ss_{\ast}}
\newcommand{\TT}{{\mathcal T}}
\newcommand{\Ta}{\TT_{\ast}}
\newcommand{\To}{\TT_{0}}
\newcommand{\Ton}{\TT_{1}}
\newcommand{\V}{{\mathcal V}}
\newcommand{\W}{{\mathcal W}}
%
%
\newcommand{\ab}{\operatorname{ab}}
\newcommand{\Aut}{\operatorname{Aut}}
\newcommand{\bg}{\operatorname{bg}}
\newcommand{\colim}{\operatorname{colim}}
\newcommand{\diag}{\operatorname{diag}}
\newcommand{\gr}{\operatorname{gr}}
\newcommand{\haut}{\operatorname{haut}}
\newcommand{\ho}{\operatorname{ho}}
\newcommand{\holim}{\operatorname{holim}}
\newcommand{\Hom}{\operatorname{Hom}}
\newcommand{\Id}{\operatorname{Id}}
\newcommand{\Image}{\operatorname{Im}\,}
\newcommand{\Ker}{\operatorname{Ker}\,}
\newcommand{\Obj}{\operatorname{Obj}\,}
\newcommand{\Prim}{\operatorname{Prim}\,}
\newcommand{\Tot}{\operatorname{Tot}}
%
%
\newcommand{\Bss}{A_{\ast\ast}}
\newcommand{\Cu}[2]{C_{{#2}}^{(#1)}}
\newcommand{\Cus}[2]{\Cu{#1}{#2}}
\newcommand{\Cuds}[1]{\Cus{#1}{\bullet}}
\newcommand{\Cds}{C_{\bullet}}
\newcommand{\Gd}{G}
\newcommand{\Jds}{J_{\bullet}}
\newcommand{\Ju}[2]{J_{{#2}}^{(#1)}}
\newcommand{\Jus}[2]{\Ju{#1}{#2}}
\newcommand{\Juds}[1]{\Jus{#1}{\bullet}}
\newcommand{\Uds}{U_{\bullet}}
\newcommand{\Vd}{V}
\newcommand{\Wd}{W}
\newcommand{\Wdd}{W_{\bullet\bullet}}
\newcommand{\Wds}{W_{\bullet}}
\newcommand{\Xd}{X}
\newcommand{\Xp}[1]{X\lra{#1}}
%
%
\newcommand{\bB}{\mathbf{B}}
\newcommand{\Di}[2]{{\EuScript D}^{#1}_{(#2)}}
\newcommand{\bK}{\mathbf{K}}
\newcommand{\bk}[1]{\mathbf{k}_{#1}}
\newcommand{\bP}[1]{\mathbf{P}_{#1}}
\newcommand{\tP}[1]{P_{#1}}
\newcommand{\Sp}{{\EuScript S}}
\newcommand{\Spp}[2]{\Sp^{#1}_{(#2)}}
%
%
\newcommand{\A}{{\EuScript A}}
\newcommand{\BA}{C\A}
\newcommand{\BW}[1]{b\WW{#1}}
\newcommand{\BG}{\BW{\Gamma}}
\newcommand{\fin}{\operatorname{fin}}
\newcommand{\init}{\operatorname{init}}
\newcommand{\vf}{v_{\fin}}
\newcommand{\vfi}[1]{v{#1}_{\fin}}
\newcommand{\vi}{v_{\init}}
\newcommand{\vin}[1]{v{#1}_{\init}}
\newcommand{\WW}[1]{W{#1}}
\newcommand{\WG}{\WW{\Gamma}}
\newcommand{\WsG}{W_{s}\Gamma}
%
%
\title{Moduli spaces of homotopy theory}
\author{David Blanc}
\address{Dept.\ of Mathematics, Univ.\ of Haifa, 31905 Haifa, Israel}
\email{blanc@math.haifa.ac.il}
\thanks{This article is dedicated to the memory of Bob Brooks,
  colleague and friend}
\date{January 25, 2004. Revised: October 3, 2004}
\begin{abstract}
The moduli spaces refered to are topological spaces whose path
components parametrize homotopy types. Such objects have been studied
in two separate contexts: rational homotopy types, in the work of
several authors in the late 1970's; and general homotopy types, in the work of
Dwyer-Kan and their collaborators. We here explain the two approaches,
and show how they may be related to each other. 
\end{abstract}
\maketitle
%
%

\sect{Introduction}
\label{cint}

The concept of ``moduli'' for a mathematical object goes back to 
Riemann, who used it to describe a set of parameters determining the
isomorphism class of a Riemann surface of a given topological type. He also
recognized that the set of all conformal equivalence classes of such
surfaces can itself be given a complex structure, although this was
only made precise by Teichm\"{u}ller, Ahlfors, Bers, and others in the
twentieth century. Similarly, the family of birational equivalence
classes of algebraic curves of each genus (over $\CC$, say) are
parametrized by an algebraic variety of moduli, as Mumford showed, 
and there are analogues in higher dimensions.

If we fix an oriented surface $S$  of genus $g$, a marked Riemann
surface is a choice of an orientation-perserving diffeomorphism of $S$
into a compact Riemann surface $X$, and the equivalence classes of
such marked surfaces form Teichm\"{u}ller space \ $\TT_{g}$, \
a complex analytic variety. The mapping class group \
$\Gamma_{g}$ \ of orientation preserving diffeomorphisms of $S$ acts
on \ $\TT_{g}$, \ and the quotient is isomorphic to the \emph{moduli space} \ 
$\M_{g}$ \ of isomorphism classes of complex structures on 
$S$ \ -- \ or equivalently, of smooth projective curves over
$\CC$. Since the action of \ $\Gamma_{g}$ \ is virtually free, and \
$\TT_{g}$ \ is contractible, \ $\M_{g}$ \ has the same rational
homology as the classifying space \ $B\Gamma_{g}$.

Note that the set of all Riemann surfaces of a given genus can be
obtained by deforming a given conformal (or equivalently, complex) 
structure on a fixed ssample surface \ $\Sigma_{g}$ \ -- \ or the 
corresponding Fuchsian group $\Gamma$. Such deformations are typically
governed by an appropriate cohomology group.

Note also that there are several levels of possible structures on a
compact surface: \ topological, differentiable, real or complex
analytic, hyperbolic, conformal, metric, and so on, which define
different types of moduli space. 
See \cite{ITaniT} or \cite{SundM} for more details and references on
the classical theory.

\subsection{Moduli spaces}
\label{sms}\stepcounter{thm}

Thus, there are a number of common themes in the moduli problems arising in
various areas of mathematics: 

\begin{enumerate}
\renewcommand{\labelenumi}{\alph{enumi})~}
\item Finding a set of parameters to describe appropriate 
equivalence classes of objects in a certain geometric
category. These may occur at varous levels, and we distinguish
between discrete parameters (such as the genus of a surface)
and the finer (continuous) \emph{moduli}. 
\item Giving the set of such equivalence classes the structure of an
object in the same category, which we may call a \emph{moduli space}.
\item The \emph{deformation} of a given structure, controlled by
cohomology, as a means of obtaining all other possible structures.
\item To pass between levels of structure, one can often quotient by 
a group action.
\end{enumerate} 

As we shall see, all these themes will find expression in the context
of homotopy theory. One aspect of the usual moduli spaces for
which we have no analogue in our setting is the important role played
by compactifications. 

\subsection{Moduli spaces in homotopy theory}
\label{smsh}\stepcounter{thm}

The homotopy theory of the classical moduli spaces of Riemann surfaces
and the corresponding mapping class groups has been studied
extensively \ -- \ see, for example, \cite{HarerSH,MWeisS,TillSC}. 
However, this is not our subject; we shall be concerned here rather
with the analogues of such moduli spaces in the category of
topological spaces, answering to the general description given above. 

The simplest examples of this approach are provided by various
mapping spaces, starting with the loop space \ $\Omega X$ \ of all
pointed maps from the circle into a topological space $X$. Here the
set \ $\pi_{0}\Omega X$ \ of components is isomorphic to the
fundamental group of $X$, and the higher homotopy groups are again
those of $X$ (re-indexed). 
More significantly, the components of space of maps from $X$ to \ $B$O \
or \ $B$U, \ the classifying spaces of the infinite orthogonal or unitary groups,
correspond to equivalence classes of (real or complex) vector bundles
over $X$, and the higher homotopy groups correspond to various
reductions of structure, by Bott periodicity.
Further examples of simple ``moduli spaces'' of topological spaces are
provided by the configuration spaces of $n$ distinct points in a
manifold $M$, which have also been studied extensively (see
\cite{FHusGT} for a comprehensive survey).

\subsection{Moduli spaces and homotopy types}
\label{smsht}\stepcounter{thm}

Our focus here will be more specifically on moduli spaces of all
homotopy types (suitably interpreted), which have been investigated from at
least two points of view over the past twenty-five years: 

\begin{enumerate}
\renewcommand{\labelenumi}{$\bullet$}
\item In work of Lemaire-Sigrist, F\'{e}lix and Halperin-Stasheff,
simply-connected rational homotopy types correspond to the 
components of a certain (infinite dimensional) algebraic variety over
$\QQ$ \ -- \ or alternatively, the quotient of such a variety by
appropriate group actions (Sec.~\ref{crhdt}).
\item The ``classification complex'' approach of Dwyer-Kan yields a 
space whose components parametrize the homotopy types of all CW complexes
(Sec.~\ref{cnms}).
\end{enumerate}

Of course, neither statement has much content, as stated; the point is
that both approaches offer (surprisingly similar) tools for
inductively analyzing the collection of all relevant homotopy types:

In each case we begin with a coarse classification of homotopy types by
algebraic invariants, such as the cohomology algebra $H$ of a space,
or the corresponding structure on homotopy groups, called a \Pa\ 
(which reduces rationally to a graded Lie algebra over $\QQ$ \ -- \ see
Sec.~\ref{cmsgla}). One might also refine this initial
classification using both homotopy and cohomology groups (Sec.~\ref{crms}). 

\begin{enumerate}
\renewcommand{\labelenumi}{(\alph{enumi})}
\item For each $H$, we have  homological classification of spaces with
this cohomology algebra. This is usually presented as an obstruction
theory, stated in terms of (algebraically defined) cohomology groups
of $H$. It can be reinterpreted as expressing the moduli space as the
limit of a tower of fibrations (Sec.~\ref{cacc}), with successive
fibers described in terms of these cohomology groups (Sec.~\ref{catt}).
\item In the rational case, we have a deformation theory approach,
which describes all homotopy types with \ $H^{\ast}(X,\QQ)\cong H$ \
by perturbing a canonical model. In fact, this was the main focus of the
rational homotopy theory approach mentioned above.
\end{enumerate}

Our main goal here is to bring out the connection between the integral
and rational cases. From the first point of view, this can be made
explicit at the algebraic level, comparing the two cohomology theories
through appropriate spectral sequences (Sec.~\ref{ccriv}).

On the other hand, the rational deformation theory has no integral
analogue, since it uses differential graded models, which exist only
for simply-connected rational homotopy types. However, it turns 
out that the deformations can be described geometrically in terms of higher 
homotopy operations (Sec.~\ref{cmhho}), and these do have integral versions. 

\begin{remark}\label{rmoduli}\stepcounter{subsection}
One could \ -- \ with some justice \ -- \  argue that the analogy between the
rational and integral ``spaces of homotopy types'' we describe here and 
classical moduli spaces is rather tenuous (although it is merely 
incidental to the point we want to make). One aspect in particular
for which we can offer no real analogue is the fact that the classical
moduli spaces have the same structure (algebraic or analytic variety,
etc.) as the objects being classified. 

In some sense, however, a homotopy type \emph{is} just the collection of
all its homotopy invariants, suitably interpreted (see \cite[I,4.9]{QuiH}).
In fact, one goal of homotopy theory is to produce a manageable
set of such invariants, sufficient for distinguishing between type. 
A starting point for such a set of invariants is the Postnikov system
of a space (and its invariants); so it is perhaps fitting that our
analysis of the space of all homotopy types is carried out by means of
its Postnikov system, and the associated collections of homotopy
invariants (higher homotopy operations, or cohomology classes) which 
appear in our obstruction theory. Thus one could choose to view these collections
of invariants themselves \ -- \ rather than the space of homotopy types \ -- \ 
as the true ``moduli object''. This is certainly the best way to understand
the connections between the rational and integral cases, and it also
might make (somewhat far-fetched) sense of the claim that the moduli
object is itself of the same kind as the things it classifies. 
We leave the philosophically-inclined reader to pursue this line of
thought at his or her discretion.
\end{remark}

\begin{ack}\label{sack}\stepcounter{subsection}
I would like to thank the referee for his or her comments, and Bill Dwyer,
Paul Goerss, Jim Stasheff, and Yves F\'{e}lix for helpful
conversations and remarks on various aspects of the subject.
\end{ack}
%
%
\sect{Rational homotopy and deformation theory}
\label{crhdt}

We first describe the deformation-theoretic approach to rational
homotopy types, starting with some basic background material:

\subsection{Rational homotopy}
\label{srht}\stepcounter{thm}

Rational homotopy theory deals with rational algebraic invariants of homotopy
types \ -- \ that is, it disregards all torsion in the homology and
homotopy groups. 
Quillen and Sullivan, respectively, proposed two main algebraic
models for the rational homotopy type of a simply-connected space $X$\vsm:

\begin{enumerate}
\renewcommand{\labelenumi}{(\roman{enumi})}
\item A differential graded-commutative algebra, or DGA \ $(A^{\ast}_{X},d)$, \ 
where \ $(A_{X}^{i})_{i=0}^{\infty}$ \ is a graded-commutative 
algebra over $\QQ$, equipped with a differential $d$ of degree \ $+1$ \ 
(which is a graded derivation with respect to the product), and \
$H^{\ast}(A_{X},d)\cong H^{\ast}(X,\QQ)$. \  Since $X$ is 
simply-connected, we may assume  \ $(A^{\ast}_{X},d)$, \ is, too \ -- \ 
that is, \ $A_{X}^{0}=\QQ$, \ $A_{X}^{1}=0$\vsm. 
\item A differential graded Lie algebra, or DGL \ $(L^{X},\partial)$, \ with \
$L^{X}=(L^{X}_{i})_{i=1}^{\infty}$ \ a positively-graded Lie algebra over
$\QQ$, and $\partial$ a graded Lie derivation of degree \ $-1$, \ such that \ 
$H_{i}(L^{X},\partial)\cong \pi_{i}\Omega X\otimes\QQ$ \ for \ $i\geq 1$, \
with the Samelson products as Lie brackets\vsm. 
\end{enumerate}

For more information see the survey \cite{FHThR}, and the original
sources \cite{QuiR,SullI}. 

\begin{defn}\label{dminmod}\stepcounter{subsection}
A DGA of the form \ $(\Lambda V,d)$, \ where \ $\Lambda V$ \ is
the free graded-commutative algebra generated by the 
graded vector space \ $V=(V^{i})_{i=0}^{\infty}$ \ is called
\emph{cofibrant}. \ A cofibrant DGA \ $(\Lambda V,d)$ \ such that \ 
$\Image(d)\subseteq\Lambda^{+}V\cdot\Lambda^{+}V$
(where \ $\Lambda^{+}V$ \ is the sub-algebra of elements in
positive degree), is called a \emph{minimal model}.
\end{defn}

A DGA or DGL map which induces an isomorphism in (co)homology \ 
is called a \emph{quasi-isomorphism}. \ 

%
%
\begin{prop}[see \protect{\cite{BLeM}}]\label{pone}\stepcounter{subsection}
Any simply-connected DGA $A$ has a minimal model with a 
quasi-isomorphism \ $q:(\Lambda V,d)\to(A,d_{A})$, \ and this is
unique up to isomorphism.  
\end{prop}

One can also define cofibrant and minimal models for DGLs, and show
that they have similar properties (cf.\ \cite[\S 21]{FHThR}).

\begin{remark}\label{rmodcat}\stepcounter{subsection}
Both the category \ $\DGAo$ \ of simply-connected DGAs, and the 
category \ $\DGLz$ \ of connected DGLs have model category
structures (see \cite{QuiH}); thus each has a concept of homotopy, and
as Quillen showed, the corresponding homotopy categories \
$\ho\DGAo$ \ and \ $\ho\DGLz$ \ are both equivalent to the rational
homotopy category \ $\ho\TT^{\QQ}_{1}$ \ of simply-connected
topological spaces (cf.\ \cite{QuiR}). The weak equivalences in \ 
$\TT^{\QQ}_{1}$ \ correspond of course to quasi-isomorphisms for DGAs or DGLs. 

Note that every simply-connected graded-commutative algebra over $\QQ$
is realizable as the rational cohomology ring of some
($1$-connected) space, and similarly, every positively graded Lie
algebra is realizable as the rational homotopy groups of a space.
\end{remark}

\subsection{The space of rational homotopy types}
\label{smrht}\stepcounter{thm}

We may therefore identify the collection \ $\M=\M^{\Ton^{\QQ}}$ \ of all rational
weak homotopy types with the set of weak homotopy types of $1$-connected
DGAs (over $\QQ$) \ -- \ i.e., equivalence classes under the
equivalence relation generated by the quasi-isomorphisms  \ -- \ and
similarly for DGLs. Moreover, there is a coarse
classification of such types, with the discrete parametrization
provided by the cohomology rings; we denote the sub-collection of $\M$
consisting of all DGAs with a given graded cohomology ring $H$ \ 
by \ $\M_{H}$. \ This has a distinguished element:

\begin{defn}\label{dformal}\stepcounter{subsection}
A DGA $A$ is called \emph{formal} if it is
quasi-isomorphic to \ $(H^{\ast}A,0)$ \ (the cohomology ring of $A$
thought of as a DGA with zero differential). Dually, a DGL $L$ is 
\emph{coformal} if it is quasi-isomorphic to \ $(H_{\ast}L,0)$.
\end{defn}

$\M_{H}$ \ was studied by Jean-Michel Lemaire and Fran{\c{c}}ois
Sigrist in \cite{LSigD}, by Yves F\'{e}lix in \cite{FelDT}, and
independently by Steve Halperin and Jim Stasheff in \cite{HStasO}. They
showed that every rational homotopy type in \ $\M_{H}$ \ can be
obtained by suitable deformations of a fixed model of the unique
formal object:

\begin{defn}\label{dbigr}\stepcounter{subsection}
The \emph{bigraded model} \ $(B,d)$ \ for $H$ is the minimal model \
$B=\Lambda Z$ \ for \ $(H,0)$ \ of Proposition \ref{pone}, \ where \
$Z=Z^{\ast}_{\ast}$, \ the graded vector space of generators, is
equipped with an additional homological grading (indicated by the
lower index).
\end{defn}

The differential takes \ $Z^{i}_{n}$ \ to \ 
$\Lambda(\oplus_{k=1}^{n-1} Z_{k})^{i+1}$, \ so one still has a
homological grading on the cohomology groups of $B$. 
The vector spaces \ $Z^{\ast}_{n}$ \ are defined (almost
canonically) by induction on $n$, so that at the $n$-th stage we kill
the cohomology in homological dimension $k$ for \ $0<k<n$. \  
$\Lambda Z_{\ast},d)$ \ is essentially the Tate-Jozefiak resolution of a 
graded commutative algebra $H$ (cf.\ \cite{TateH,JozeT}) \ -- \
i.e., a minimal ``cellular'' resolution. See \cite[\S 3]{HStasO} for the details. 

The deformation consists of a perturbation of the 
differential $d$ to \ $D=d+D'$ \ -- \ of course, this is done in such
a way that the cohomology is unchanged, and  \ $D\circ D=0$. \ 
We can no longer expect $D$ to respect both gradings, but if we think
of $B$ as being only \emph{filtered} by \
$\F_{n}B:=\oplus_{k-0}^{n}~(\Lambda Z^{\ast})_{k}$, \ then the
perturbation is such that \ $D'$ \ takes \ $Z_{n}$ \ into \ $\F_{n-2}B$. \ 

%
%
\begin{thm}[cf.\ \protect{\cite[Theorem 4.4]{HStasO}}]
\label{tone}\stepcounter{subsection}
If \ $(A,d_{A})$ \ is any DGA with \ $H^{\ast}(A,d_{A})\cong H$, \ one can
deform the bigraded model \ $(B,d)$ \ for \ $(H,0)$ \ into a filtered
model \ $(B,D)$ \ quasi-isomorphic to \ $(A,d_{A})$, \ and
this is unique up to isomorphism.  
\end{thm}

\subsection{The variety of deformations}
\label{svdef}\stepcounter{thm}

The deformations $D$ of $d$ are determined by successive choices
of linear transformations \ $D'_{n}:Z_{n}\to\F_{n-2}\Lambda Z_{<n}^{\ast}$, \ 
satisfying a quadratic algebraic condition (to ensure that \ $D^{2}=0$). \ 
Thus the collection of all possible perturbations of the given
bigraded model \ $(B,d)$ \ constitute an (infinite dimensional) 
algebraic variety over $\QQ$, \ which we denote by \ $\V_{H}$. \ 
Every simply-connected rational homotopy type is represented by at
least one point of \ $\V_{H}$; \ each such point actually represents an
``augmented'' DGA $X$, equipped with a specific isomorphism \ 
$\phi:H^{\ast}(X;\QQ)\to H$.

However, one can have more than one point corresponding to each augmented
rational homotopy type: the points of \ $\M_{H}$ \ 
correspond bijectively with the orbits of \ $V_{H}$ \ under the action
of an appropriate group, as in
\cite[Thm.~4]{LSigD} or \cite[Ch.~1]{FelDT} \ (note that these two
superificially different descriptions are in fact the same, by \cite{FelM}):

In the approach described in \cite{SStasD}, we first divide by the
action of a certain pro-unipotent algebraic group \ $G_{H}$ \ to obtain the
set \ $V_{H}/G_{H}$ \ of augmented rational homotopy types as above.
We must then further mod out by \ $\Aut(H)$ \ to get the plain
(i.e., unaugmented) rational homotopy types, \ so that \ 
$\M_{H}\cong\Aut(H)\backslash V_{H}/G_{H}$. \  
In \cite[\S 6]{SStasD} it is then shown that the quotient set \
$V_{H}/G_{H}$ \ can also be described as the set of path components of the
Cartan-Chevalley-Eilenberg ``standard construction'' \ $C(L)$ \ 
(cf.\ \S \ref{sqe}) on a certain Lie algebra $L$ (namely, the algebra
of derivations of the bigraded model for $H$).

Note that though \ $\M_{H}$ \ can be uncountable, it can also be small: a 
DGA \ $(A,d)$ \ is called \emph{intrinsically formal} if \ 
$\M_{H^{\ast}A}$ \ consists of a single point.

\subsection{Obstruction theory}
\label{sot}\stepcounter{thm}

We have seen that the collection \ $\M^{\Ton^{\QQ}}$ \ of all
simply-connected rational homotopy types indeed exhibits the basic
characteristics of a moduli space listed in \S \ref{sms}. However, to
make this useful, we need further information about its overall structure.  

This is provided by Halperin and Stasheff in the form of an
obstruction theory for realizing a given isomorphism of rational cohomology
algebras \ $\phi:H^{\ast}(A,d_{A})\to H^{\ast}(B,d_{B})$ \ by a DGA
map (or equivalently, by a map of the corresponding rational spaces). 
Since such a map induces a homotopy equivalence between the respective 
cofibrant models, this can be thought of as a method for
distinguishing between the homotopy types in \ $\M_{H}$.

The obstruction theory, which only works in full generality if the
DGAs in question are of finite type, was presented in \cite[\S 5]{HStasO}
in terms of a sequence of elements \ 
$\OO_{n}(f)\in\Hom(Z_{n+1},H(B))$ \ -- \ or rather, in a quotient of
this group. \ In \cite[\S 4]{StaRH}, these are re-interpreted as
elements in \ $\bigoplus_{p}~H^{p+1}(B;\pi_{p}(H(A)))$.

\begin{remark}\label{rcoh}\stepcounter{subsection}
By analogy with classical deformation theory, we would
expect the deformations of the bigraded model to be ``controlled'' by a
suitable differential graded Lie algebra \ $(L,d)$, \ or perhaps by
its cohomology. In fact, F\'{e}lix observed that the Halperin-Stasheff
obstructions \ $\OO_{n}(f)$ \ can be thought of as lying in the
F-cohomology group \ $_{F}H^{1}_{n}(B,d)$ \ of the bigraded
model \ $(B,d)$ \ for \ $(H,0)$ \ (see \cite[Ch.~5, Prop.~7]{FelDT}). \ 
Moreover, these groups can be identified with the tangent (or Harrison)
cohomology of the graded algebra $H$ (cf.\ \cite[Ch.~3, Prop.~4]{FelDT}), 
which is the same as the Andr\'{e}-Quillen cohomology, and this can
indeed be calculated as the cohomology of a suitable DGL (cf.\
\cite[\S 4]{SStasL}).

Moreover, both Schlessinger and Stasheff (see \cite[\S 5]{StaRH}) and 
F\'{e}lix (cf.\ \cite[Ch.~4-5]{FelDT}) obtain results connecting 
formality and intrinsic formality of a DGA with its tangent
cohomology, though the relation does not work both ways (see
\cite[Ch.~4, Ex.~2]{FelDT} ).

We shall give a uniform treatment of this cohomological description in
a more general context in Section \ref{catt}. However, for technical
reasons (inter alia, convergence of second- versus first-quadrant 
spectral sequences), this is done for the Eckmann-Hilton dual \ -- \ that is,
in terms of \Pa s (the integral version of graded Lie algebras over
$\QQ$). The mod $p$ version \ -- \ that is, an obstruction theory for
realizing unstable (co)algebras over the Steenod algebra \ -- \
appears in \cite{BlaR}.
\end{remark}

%
%
\sect{Moduli spaces of graded Lie algebras}
\label{cmsgla}

The description of the moduli space of rational homotopy types in terms of DGLs,
rather than DGAs, has received less attention. However, as noted above,
it has certain technical advantages, especially in the integral case
(where we need no longer restrict to simply-connected spaces, incidentally). 
The dualization to DGLs is quite straightforward, but we shall
describe it here from a slightly different point of view. 

\subsection{Deformations of DGLs}
\label{sddgl}\stepcounter{thm}

The deformation theory for DGLs proceeds precisely along the lines
indicated in \S \ref{smrht}\textit{ff.}: \ given a positively-graded
Lie algebra \ $P=P_{\ast}$ \ over $\QQ$, consider the collection \
$\M_{P}$ \ of all rational homotopy types of simply-connected spaces
$X$  \ -- \ or equivalently, of DGLs \ $(L,\partial)$ \ -- \ such that \ 
$\pis\Omega X\otimes\QQ=H_{\ast}(L,\partial)\cong P$ \ (as Lie algebras). 

Again there is a distinguished coformal element  \ $(P,0)$ \ in \ $\M_{P}$ \ 
(Def.~\ref{dformal}). Its DGL minimal model \ $(B,\partial)$ \
can be given an additional homological grading, and
the differential in this bigraded model can be deformed so as to
obtain a filtered model for any $X$ in \ $\M_{P}$. \ See the sketch in
\cite{HaralP}; the details appear in the unpublished thesis \cite{OukiH}.

\subsection{The resolution model category}
\label{srmc}\stepcounter{thm}

Since the category \ $\DGLz$ \ is not abelian, from the point of view
of homotopical algebra it is more natural to consider simplicial 
resolutions of objects in this category, rather than ``chain complex''
resolutions (which is what the homological grading of the bigraded
model gives us). Moreover, this generalizes more naturally to the
integral setting, where differential graded models are no longer
available (although see \cite{MandE} for a way around this).

Recall that a \emph{simplicial object} \ $X_{\bullet}$ \ over a category $\C$
is a functor \ $X:\Delta^{op}\to\C$, \ where \ $\Delta$ \ is the
category of finite sequences \ $\bn=(0,1,\dotsc, n)$ \ ($n\in\NN$), \
with order-preserving maps as morphisms. The category of simplicial
objects over $\C$ will be denoted by \ $s\C$. 

Of course, \ $\C=\DGLz$ \ (as well as \ $\DGAo$) \ is itself a model
category (as shown by Quillen in \cite{QuiR}), so in considering
simplicial DGLs we can take advantage of the procedure described by
Dwyer, Kan, and Stover for generating a model category structure on \
$s\C$. \ Essentially, this consists of: 

\begin{enumerate}
\renewcommand{\labelenumi}{$\bullet$}
\item choosing a certain set \ $\{M_{\alpha}\}_{\alpha\in A}$ \ of
``good'' objects in $\C$ (corresponding to the projectives of
homological algebra); 
\item creating simplicial objects \ $\Sigma^{n}M_{\alpha}$ \ for each \ 
$n\in\NN$ \ and \ $\alpha\in A$ \ by placing \ $M_{\alpha}$ \ 
in simplicial dimension $n$; \ and using these to
\item define weak equivalences (which are detected by maps out of \
$\Sigma^{n}M_{\alpha}$), \ and 
\item construct the cofibrant objects (corresponding to 
projective resolutions), and more generally all cofibrations.
\end{enumerate}

See \cite{DKStE} for more details; as extended by Bousfield (who calls it a
``resolution model category'') in \cite{BousCR}, this procedure works
in great generality. In particular, it includes the model categories
for universal algebras defined by Quillen in \cite[II, \S 4]{QuiH}, as
well as the graded version (see, e.g., \cite[\S 2.2]{BStoG}). The
original example considered by Dwyer, Kan, and Stover was that of
simplicial spaces; this can be readily modified to cover simplicial
rational spaces (equivalently, simplicial DGLs or DGAs).

\subsection{Simplicial DGLs}
\label{ssdgl}\stepcounter{thm}

In applying the above procedure to \ $\C=\DGLz$, \ we have a
natural candidate for the ``good objects'' \ -- \ namely, the DGL spheres \ 
$\{\Sp^{n}\}_{n=2}^{\infty}$ \  (which are just minimal models for \ 
$\{S^{n}_{\QQ}\}_{n=2}^{\infty}$ \ -- \ i.e., \ $\Sp^{n}=\Spp{n}{x}$ \ 
is a free graded Lie algebra on a single generator $x$ in dimension $n$).

This determines what we mean by a (simplicial) \emph{resolution} of a
given DGL \ $(L,\partial)$ \ -- \ namely, a cofibrant 
simplicial DGL \ $\Vd$, \ with each \ $V_{n}$ \ homotopy equivalent to a
coproduct \ $\coprod_{i\in I_{n}}~\Sp^{n_{i}}$, \ such that the 
augmented simplicial DGL \ $\Vd\to L$ \ is weakly contractible. 
A map \ $\varphi:\Vd\to\Xd$ \ of simplicial DGLs is a weak equivalence
if and only if after taking (internal) homology in each simplicial
dimension, the resulting map of graded simplicial abelian groups \ 
$\varphi_{\ast}:H_{\ast}\Vd\to H_{\ast}\Xd$ \ induces an isomorphism \ 
$H_{n}(\varphi_{\ast}):H_{n}(H_{\ast}(\Vd))\xra{\cong} H_{n}(H_{\ast}(\Xd))$ \ 
for each \ $n\geq 0$.

By \cite[Thm.~2.1]{QuiR}, we know that the model category \ $\DGLz$ \
is equivalent to the category \ $s\Lie_{1}$ \ of connected simplicial
Lie algebras; similarly for bigraded differential Lie algebras and
simplicial graded Lie algebras by \cite[Prop.\ 2.9]{BlaHR}. In fact,
it turns out that there is also a one-to-one correspondence between
filtered DGLs and (cofibrant) simplicial DGL resolutions (see
\cite[Prop.\ 5.13]{BlaHR}), so the DGL deformation theory mentioned in
\S \ref{sddgl} can be translated into the language of simplicial
resolutions. 

\subsection{Moduli and realizations}
\label{smr}\stepcounter{thm}

One can analyze \ $\M_{P}$ \ in terms of a \emph{realization problem}:
which distinct rational homotopy types of $1$-connected topological
spaces have the given graded Lie algebra $P$ as their rational 
homotopy groups?  In view of the equivalences of categories noted
above, we can replace the DGL bigraded model \ $(B,\partial)$ \ for \
$(P,0)$ \ by a free simplicial resolution \ $\Gd\to P$ \ of graded Lie
algebras, and observe that if $P$ itself can be realized by a rational
space \ $X$ \ (and so by a DGL \ $(L^{X},\partial)$), \ then \ $\Gd$ \
can be realized by a \emph{simplicial} space (or DGL) \ $\Vd$. \
Moreover, any weak equivalence \ $\varphi:\Vd\to\Wd$ \ of simplicial
spaces (or DGLs) induces a weak equivalence \ 
$\varphi_{\#}:\pis\Vd\otimes\QQ\to\pi_{\ast}\Wd\otimes\QQ$, \
and since \ $\pis\Vd\otimes\QQ=\Gd$ \ is a resolution of $P$, \
this is just a choice of an isomorphism between two isomorphic copies of $P$. 

As we shall see below (\S \ref{cnms}), there is a ``moduli space'' associated
to any model category, and the various levels of classification which
are a feature of moduli spaces in general (\S \ref{sms}(i)) can be
understood in terms of ``structure-reducing'' functors between the
corresponding model categories \ -- \ here exemplified by \
$\pis(~)\otimes\QQ:\TT_{1}\to\gr\Lie$ \ (extended to simplicial objects).

%
%

\sect{Moduli and higher homotopy operations}
\label{cmhho}

Ideally, we want not only an algebraic description of the components of \
$\M_{H}$ \ in terms of the cohomology of the graded algebra $H$, say,
as in \S \ref{sot}, but also a geometric interpretation of this
description \ -- \ with the potential for lifting to the integral case. 

Note that there is an obvious candidate for geometric invariants to 
distinguish between inequivalent rational spaces with the same
cohomology \ -- \ namely, Massey products, and their higher order
analogues (see \cite{MassN,KraiM}).
In fact, Halperin and Stasheff's original motivation for their
obstruction theory was to make precise a folk theorem stating that a space
is formal if and only if all higher order Massey products \ -- \ 
or more generally, all higher-order cohomology operations \ -- \ vanish.

One difficulty is in the very definition of higher cohomology
operations (see, e.g., \cite{MaunC}). This can be bypassed to some
extent by re-interpreting them in terms of differentials in a spectral
sequence (cf.\ \cite[\S 7]{HStasO}) \ -- \  at the price of losing both
their geometric content, and the ability to identify the same higher
operation in different contexts. 

Because their Eckmann-Hilton dual \ -- \ namely, higher homotopy
operations \ -- \ is somewhat more intuitive than the cohomology
version, we shall concentrate on these. A precise definition requires
some care; in the following sections, we present the version of \cite{BMarkH}: 

\begin{defn}\label{dlatt}\stepcounter{subsection}  
A \emph{lattice} is a finite directed non-unital category $\Gamma$
(that is, we omit the identity maps), equipped with two objects \
$\vi=\vi(\Gamma)$ \ and \ $\vf=\vf(\Gamma)$ \ with a unique map \
$\phi_{\max}:\vi\to\vf$, \ and for every \ $w\in V\DEF\Obj\Gamma$, \ there is 
at least one map from \ $\vi$ \ to $w$, and at least one from $w$ to \ $\vf$. 
A composable sequence of $k$ arrows in $\Gamma$ will be called a
$k$-\emph{chain}.
\end{defn}

\subsection{The $W$-construction}
\label{swc}\stepcounter{thm}

Given a lattice $\Gamma$, one can define a new category \ $\WG$ \ 
enriched over cubical sets, with the same set of objects $V$ (cf.\
\cite[III, \S 1]{BVogHI}):

For \ $u,v \in V$, \ let \ $\Gamma_{n+1}(u,v)$ \ be the set of \ 
$(n+1)$-chains from $u$ to $v$ in $\Gamma$, and \ 
$\WG(u,v) := \bigsqcup_{n \geq 0} \Gamma_{n+1}(u,v) \times I^{n}/\sim$, \ 
where $I$ is the unit interval. Write \ 
$f_{1}\ccsub{t_{1}}f_{2}\dotsb f_{n}\ccsub{t_{n}}f_{n+1}$ \ for \ 
$\lra{u\xra{f_{n+1}} v_{n}\dotsb\to v_{1}\xra{f_{1}}v}
\times (t_{n},\dotsc,t_{1})$ \ in \ $\Gamma_{n+1}(u,v) \times I^{n}$; \ 
then the relation $\sim$ is generated by
$$
f_{1}\ccsub{t_{1}}f_{2}\dotsb f_{n}\ccsub{t_{n}}f_{n+1}~\sim~
f_{1} \circ_{t_1} \dotsb \circ_{t_{i-1}} 
(f_i f_{i+1}) \circ_{t_{i+1}} \dotsb \circ_{t_n} \circ f_{n+1}\hsp
\text{if}\hsm t_{i} =0
$$
\noindent for \ $1 \leq i \leq n$, \ where \ 
$(f_{i} f_{i+1})$ \ denotes \ $f_{i}$ \ composed with \ $f_{i+1}$.

The categorial composition in \ $\WG$ \ is given by the concatenation:
$$
(f_{1}\circ_{t_{1}} \dotsb \circ_{t_{l}} f_{l+1}) \circ 
(g_{1}\circ_{u_{1}} \dotsb \circ_{u_{k}} g_{k+1})\DEF 
(f_{1}\circ_{t_{1}} \dotsb \circ_{t_{l}} f_{l+1} \circ_{1} g_{1}
\circ_{u_{1}} \dotsb \circ_{u_{k}} g_{k+1}). 
$$
We write \ $\PP\DEF\WG(\vi,vf)$ .

\begin{defn}\label{dbas}\stepcounter{subsection}  
The \emph{basis} category \ $\BG$ \ for a lattice $\Gamma$ is defined
to be the cubical subcategory of \ $\WG$ \ with the same objects, and
with morphisms given by \ $\BG(u,v)\DEF\WG(u,v)$ \ if \
$(u,v)\neq(\vi,\vf)$, \ while 
$$
\BG(\vi,\vf):=
\bigcup \left\{\alpha\circ\beta\ | \ \beta\in\WG(\vi,w),\
\alpha\in\WG(w,\vf),\ \vi\neq w\neq\vf \right\},
$$
so that \ $b\PP\DEF\BG(\vi,\vf)$ \ consists of all decomposable morphisms.
\end{defn}

%
%
\begin{fact}[\protect{\cite[Proposition 2.15]{BMarkH}}]
\label{fone}\stepcounter{subsection}  
For any lattice $\Gamma$, \ $\WG(\vi,\vf)$ \ is isomorphic to the cone
on its basis, with vertex corresponding to the unique maximal $1$-chain.
\end{fact}

We can use the $W$-construction to define a higher homotopy operations,
as follows:

\begin{defn}\label{dhho}\stepcounter{subsection}  
Given \ $\A:\Gamma\to\ho\Ta$ \ for a lattice $\Gamma$ as above, 
the corresponding \emph{higher order homotopy operation} \ is the subset \ 
$\llrr{\A}\subset [\A (\vi) \rtimes b\PP,\A (\vf)]_{\ho\Ta}$ \ 
of the homotopy equivalence classes of maps  \ 
$$
\BA|_{b\WG(\vi,\vf)}: b\WG(\vi,\vf) = 
b\PP \longrightarrow\Ta(\A(\vi),\A(\vf))
$$ 
induced by all possible continuous functors \ $\BA:\BG\to\Ta$ \ such that \ 
$\pi\circ\BA  = \A\circ(\varepsilon\rest{\BG})$, \ 
where the half-smash is defined \ 
$X\rtimes K\DEF(X\times K)/(\{\ast\}\times K)=X\wedge K_{+}$. \ 

$\llrr{\A}$ \ is said to \emph{vanish} if it contains the homotopy
class of a constant map \ $b\PP\longrightarrow\Ta(\A(\vi),\A(\vf))$.
\end{defn}

Given a lattice $\Gamma$, a \emph{rectification} of a functor \
$\A:\Gamma\to\ho\Ta$ \ is a functor \ $F:\Gamma\to\Ta$ \ with \
$\pi\circ F$ \ naturally isomorphic to $\A$ (where \
$\pi:\Ta\to\ho\Ta$ \ is the obvious projection functor).  

%
%
\begin{prop}[\protect{\cite[Thm.~3.8]{BMarkH}}]
\label{ftwo}\stepcounter{subsection}  
$\A$ has a rectification if and only if the homotopy operation \
$\llrr{\A}$ \ vanishes (and in particular, is defined) 
\end{prop}

With this at hand, we can now reformulate the deformation
theory for the bigraded DGL model for \ $(P,0)$ \ in the following form:

%
%
\begin{thm}[\protect{\cite[Thm.~7.14]{BlaHR}}]
\label{ttwo}\stepcounter{subsection}
For each connected graded Lie algebra $P$ of finite type, there is a
tree \ $T_{P}$, \ with each node $\alpha$ indexed by a cofibrant DGL \ 
$L^{\alpha}$ \ such that \ $H_{\ast}L^{\alpha}\cong P$, \ 
starting with \ $L^{0}\simeq(P,0)$ \ at the root $0$. \ For each
node $\alpha$ there is an integer \ $n_{\alpha}>0$ \ such 
that \ $L^{\alpha}$ \ agrees with its successors in degrees 
$\leq n_{\alpha}$, \ with \ $n_{0}=1$, \ so that the sequential 
colimit \ $L^{(\infty)}=\colim_{k}L^{(k)}$ \ along any branch is well
defined, and for any rational homotopy type in \ $\M_{P}$ \ there
exists such a tree. 
Furthermore, for each two immediate successors \ $\beta,\gamma$ \ of a node 
$\alpha$, there is a higher homotopy operation \ 
$\llrr{\beta,\gamma}\subseteq [\Sp^{n},L]$ \ which vanishes 
for \ $L=L^{\beta}$, \ but not for \ $L=L^{\gamma}$ \ (or conversely).
\end{thm}

Thus we have reformulated the deformation of the bigraded model as an obstruction
theory for distinguishing between the various realizations of a given
rational graded Lie algebra $P$ in terms of higher homotopy
operations. This can also be done integrally (cf.\ \cite{BlaHH}),
and one could identify these operations with appropriate cohomology
classes in the dual version of the Halperin-Stasheff obstructions (see 
\cite{BlaAI}).

\subsection{Example of a DGL deformation}
\label{edgld}\stepcounter{thm}

We shall not explain how the higher operations \ $\llrr{\beta,\gamma}$ \ 
are defined, since the construction is rather technical; see \cite[\S
4]{BlaHR} for the details. Instead, the reader might find the
following example instructive\vsm :

Consider the graded Lie algebra \ 
$P=\bL{a_{1},b_{1},c_{2}}/I$, \ where \ $\bL{T_{\ast}}$ \ is
the free graded Lie algebra on a graded set of generators \
$T_{\ast}$, \ (with the subscripts indicating the dimension), and $I$
is here the Lie ideal generated by \ $[a,a]$ \ and \ $[[c,a],[b,a]]$\vs.

%
\noindent\textbf{Step I:} \ \ \ 
The minimal model for the coformal DGL \ $(P,0)\in\DGLz$ \ is \ 
$(B,\partial_{B})$, \ where $B$ in dimensions $\leq 7$ \
is \ $\bL{a_{1},b_{1},c_{2},x_{3},y_{5},w_{6},z_{7}}$, \ with \ 
$\partial_{B}(x)=[a,a]$, \ $\partial_{B}(y)=3[x,a]$, \ 
$\partial_{B}(w)=[[c,a],[b,a]]$, \ and \ 
$\partial_{B}(z)=4[y,a]+3[x,x]$. \ 
The bigraded model \ $\Bss$ \ is obtained from $B$ by
introducing an additional (homological) grading: \ $a,b\in A_{0,1}$, \
$c\in A_{0,2}$, \ $x\in A_{1,3}$, \ $w\in A_{1,6}$, \ $y\in A_{2,5}$, \
$z\in A_{3,7}$, \ and so on\vs .

%
\noindent\textbf{Step II:} \ \ \ 
The free simplicial DGL resolution \ $\Cds\to L=(\LL,0)$ \ 
may be described in homological dimensions $\leq 3$ \ as follows\vsm :

\begin{enumerate}
\renewcommand{\labelenumi}{(\arabic{enumi})}
\item $\Cus{0}{0}$ \ is the DGL coproduct of \ 
$\Spp{1}{a^{(0)}}=\Spp{1}{\lra{a}}$, \ 
$\Spp{1}{b^{(0)}}=\Spp{1}{\lra{b}}$, \ and \ 
$\Spp{2}{c^{(0)}}=\Spp{2}{\lra{c}}$.
\item $\Cus{0}{1}=\Spp{2}{x^{(0)}}\amalg \Spp{6}{w^{(0)}}$, \ where \ 
$x^{(0)}=\lra{[\lra{a},\lra{a}]}$ \ and \ 
$w^{(0)} = \lra{[[\lra{c},\lra{a}],[\lra{b},\lra{a}]]}$.
\item $\Cus{0}{2}$ \ consists of \ $\Spp{5}{y^{(0)}}$, \ where \ 
$y^{(0)} = \lra{3[\lra{[\lra{a},\lra{a}]},\lra{\lra{a}}]}$.
\item $\Cus{0}{3}$ \ consists of \ 
$\Spp{6}{z^{(0)}}$, \ where \ 

$z^{(0)}=
\lra{4[\lra{3[\lra{[\lra{a},\lra{a}]}, \lra{\lra{a}}]},\lra{\lra{\lra{a}}}] + 
6[\lra{[\lra{\lra{a}},\lra{\lra{a}}]},\lra{\lra{[\lra{a},\lra{a}]}}]}$\vsm.
\end{enumerate}

In analogy with the DGL sphere of \S \ref{ssdgl}, the DGL $n$-disk \ 
$\Di{n}{x}$ \ is the free graded Lie algebra on two generators: $x$ in
dimension $n$ and \ $\partial(x)$ \ in dimension $n-1$.
With this notation, for \ $\Cuds{1}$ \ we need in addition\vsm:

\begin{enumerate}
\renewcommand{\labelenumi}{(\arabic{enumi})}
\item $\Di{3}{x^{(1)}}\hra\Cus{1}{0}$ \ with \
$\partial_{W}(x^{(1)})=d_{1}(x^{(0)})=\lra{[a,a]}$. 
\item $\Di{6}{y^{(1)}}\hra\Cus{1}{1}$ \ with \ 
$y^{(1)}=\lra{3[\lra{x},\lra{a}]}$ \ and \  
$\partial_{W}(y^{(1)})=d_{2}(y^{(0)})=\lra{3[\lra{[a,a]},\lra{a}]}$. 
\item $\Di{7}{z^{(1)}}\hra\Cus{1}{2}$ \ with \
$z^{(1)} = \lra{4[\lra{3[\lra{x},\lra{a}]},\lra{\lra{a}}] + 
6 [\lra{[\lra{a},\lra{a}]},\lra{\lra{x}}]}$ \ and \ 
$\partial_{W}(z^{(1)}) = d_{3}(z^{(0)})=
\lra{4[\lra{3[\lra{[a,a]}, \lra{a}]},\lra{\lra{a}}] + 
6 [\lra{[\lra{a},\lra{a}]},\lra{\lra{[a,a]}}]}$\vsm .
\end{enumerate}

For \ $\Cuds{2}$ \ we need in addition \ 

\begin{enumerate}
\renewcommand{\labelenumi}{(\arabic{enumi})}
\item $\Di{7}{y^{(2)}}\hra\Cus{2}{0}$ \ with \
$\partial_{W}(y^{(2)})=d_{1}(y^{(1)})=\lra{3[x,a]}$. 
\item $\Di{8}{z^{(2)}}\hra\Cus{2}{1}$ \ with \ 
$z^{(2)} = \lra{4[\lra{y},\lra{a}] + 3 [\lra{x},\lra{x}]}$ \ and \ 
$\partial_{W}(z^{(2)})=d_{2}(z^{(1)})=$

$\lra{4[\lra{3[x,a]},\lra{a}]+6[\lra{[a,a]},\lra{x}]}$\vsm .
\end{enumerate}

For \ $\Cuds{3}$ \ we must add \ $\Di{9}{z^{(3)}}\hra\Cus{3}{0}$ \ 
with \ $\partial_{W}(z^{(1)})=d_{1}(z^{(2)})= \lra{4[y,a]+3[x,x]}$\vs.

%
\noindent\textbf{Step III:} \ \ \ 
So far, we have constructed the DGL simplicial resolution of the
coformal DGL \ $(B,\partial_{B})\simeq (P,0)$. \ Now we consider a
different DGL having the same homology (i.e., a non-weakly equivalent 
topological space with the same rational homotopy groups)\vsm :

Consider the DGL \ $(B',\partial_{B'})$ \ where \ 
$B':=\bL{a_{1},b_{1},c_{2},x_{3},y_{5},z_{7},\dotsc}$, \ 
is a free Lie algebra with \  
$\partial_{B'}(x)=[a,a]$, \ $\partial_{B'}(y)=3[x,a]-[[b,a],c]$, \ 
$\partial_{B'}(z)=4[y,a]+3[x,x]$, \ and so on.

Here \ $P\cong H_{\ast}(B)=\bL{a_{1},b_{1},c_{2}}/\lra{[a,a],[[c,a],[b,a]]}$, \ 
so the bigraded model for \ $(P,0)$ \ is \ $(\Bss,\partial_{B})$ \ 
of Step II above, \ and the filtered model is obtained from it by setting \ 
$D(y)=3[x,a]-[[b,a],c]$ \ and \ $D(z)=4[y,a]+3[x,x]-4w+2[[x,b],c]$\vsm .

The corresponding free simplicial DGL resolution is obtained from \ $\Cds$ \ of
Step II by making the following changes\vsm :

\begin{enumerate}
\renewcommand{\labelenumi}{(\arabic{enumi})}
\item Set \ 
$y^{(1)}\DEF\lra{3[\lra{x},\lra{a}]-[[\lra{b},\lra{a}],\lra{c}]}$, \ with \ 
$$
\partial_{W}(y^{(1)})=d_{2}(y^{(0)})=
\lra{3[\lra{[a,a]},\lra{a}]}
$$
as before \ (but now \ $\partial_{W}(y^{(2)})=\lra{3[x,a]-[[b,a],c]}$, \ 
of course). 
\item Set \ 
\begin{equation*}
\begin{split}
z^{(1)}~\DEF~& \langle 4[\lra{3[\lra{x},\lra{a}]},\lra{\lra{a}}] + 
6 [\lra{[\lra{a},\lra{a}]},\lra{\lra{x}}] \\
&- 4 \lra{[[\lra{c},\lra{a}],[\lra{b},\lra{a}]]} +  
2 [[\lra{[\lra{a},\lra{a}]},\lra{\lra{b}}],\lra{\lra{c}}]\rangle
\end{split}
\end{equation*}
(with \ $\partial_{W}(z^{(1)})$ \ unchanged).
\item Set \ 
$z^{(2)}\DEF\langle 4[\lra{y},\lra{a}] + 6 [\lra{x},\lra{x}] - 4 \lra{w} +
2 [[\lra{x},\lra{b}],\lra{c}]\rangle $, \ with \
\begin{equation*}
\begin{split}
\partial_{W}(z^{(2)})~=~&d_{3}(z^{(1)})~=~
\langle 4[\lra{3[x,a]-[[b,a],c]},\lra{a}] \\
&+ 6 [\lra{[a,a]},\lra{x}] -  
4 \lra{[[c,a],[b,a]]} + 2 [[\lra{[a,a]},\lra{b}],\lra{c}] \rangle
\end{split}
\end{equation*}
\item Finally, \ $\partial_{W}(z^{(3)})=\lra{4[y,a]+3[x,x]-4w+2[[x,b],c]}$\vs.
\end{enumerate}

%
\noindent\textbf{Step IV:} \ \ \ In order to define the higher
homotopy operations which distinguish \ $B'$ \ from $B$, observe that
the required half-smashed are defined as follows\vsm :

$\Spp{2}{x}\rtimes D[1]=(\bL{X_{\ast}},\partial')$, \ where \ 
$X_{2}=\{(x,(d_{0})),(x,(d_{1}))\}$, \ $X_{3}=\{(x,(\Id))\}$, \ 
and \ $\partial'(x,(\Id))=(x,(d_{0}))-(x,(d_{0}))$.

Similarly, \ $\Spp{3}{y}\rtimes D[2]=(\bL{Y_{\ast}},\partial')$, \ where \ 
$Y_{3}=\{(y,(d_{0}d_{1})),(y,(d_{0}d_{2})),(y,(d_{1}d_{2}))\}$, \ 
$Y_{4}=\{(y,(d_{0})),(y,(d_{1})),(y,(d_{2}))\}$, \ and \ 
$Y_{5}=\{(y,(\Id))\}$, \ 
with \ $\partial'(y,(\Id))=-(y,(d_{0}))+(y,(d_{1}))-(y,(d_{2}))$, \ 
$\partial'(y,(d_{0}))=-(y,(d_{0}d_{1}))+(y,(d_{0}d_{2}))$, \ 
$\partial'(y,(d_{1}))=-(y,(d_{0}d_{1}))+(y,(d_{1}d_{2}))$, \ and \ 
$\partial'(y,(d_{2}))=-(y,(d_{0}d_{2}))+(y,(d_{1}d_{2}))$\vs.

%
\noindent\textbf{Step V:} \ \ \ 
For the DGL \ $B'$ of Step III, with \ $\Cds$ \ as in Step II, \ we define \
$h_{0}=f_{0}:C_{0,\ast}\to B'$ \ by setting \ $f_{0}(\lra{a})=a$, \ 
$f_{0}(\lra{b})=b$, \ $f_{0}(\lra{c})=c$, \ 
and \ $f_{0}=0$ \ for all other disks in \ $C_{0,\ast}$ \ (so for
example \ $f_{0}(\lra{[a,a]})=0$).

Thus on \ $\Spp{2}{x^{(0)}}\rtimes D[1]$ \ we have \ 
$h_{1}(x^{(0)},(d_{0}))=[a,a]$ \ and \ $h_{1}(x^{(0)},(d_{1}))=0$, \ 
(see \cite[Def.~4.7]{BlaHR}), so we must choose \
$h_{1}(x^{(0)},(\Id))=x\in B'_{3}$. 

Now on \ $\Spp{3}{y^{(0)}}\rtimes D[2]$ \ we have \ 
$h_{2}(y^{(0)},(d_{0}d_{1}))= h_{2}(y^{(0)},(d_{0}d_{2}))=
h_{2}(y^{(0)},(d_{1}d_{2}))=0$, \ and \ 
$h_{2}(y^{(0)},(d_{0}))=(h_{1}\circ d_{0})(y^{(0)},(d_{0})=3[x,a]$, \ 
while \ $h_{2}(y^{(0)},(d_{1}))=h_{2}(y^{(0)},(d_{2}))=0$\vsm .

Thus in Definition \ref{dhho} we shall be interested in a lattice
$\Gamma$ corresponding to the face maps of a simplicial object, with \ 
$\BG\cong\partial D[2]$. \ We have just defined a continuous functor \ 
$\BA:\BG\to\Ta$ \ for \ $\Cds$, \ as required there, 
and the resulting secondary operation is \
$\llrr{2,y^{(0)}}=\{\lra{3[x,a]}]\}\subseteq H_{4}B'$ \ (remember that
the homology of a DGL corresponds to the rational homotopy groups of
the associated space). \ Since \ $3[x,a]$ \ does not bound in \ $B'$, \ 
$\llrr{2}$ \ does not vanish, and we have found an obstruction to the
coformality of \ $B'$\vsm. 

Unfortunately, DGL higher homotopy operations such as \
$\llrr{2,y^{(0)}}$ \ are unsatisfactory, in as much as one cannot
translate them canonically into integral homotopy operations. 
One way to avoid this difficulty is to use more general non-associative algebras,
rather than Lie algebras, as our basic models. 
Thus, consider the category  \ $\DGN$ \ of non-associative
differential graded algebras. A DGN whose homology happens
to be a graded Lie algebra will be called a \emph{Jacobi}
algebra. Every DGL has a free Jacobi model, and these can be
used to resolve any DGL in \ $s\DGN$ \ (see \cite[\S 7]{BlaHR})\vs.

%
\noindent\textbf{Step VI:} \ \ \ Consider a Jacobi model $G$ for the
coformal DGL $P$ of Step I. We construct the minimal Jacobi resolution \ 
$\Jds\to G$ \ corresponding to the bigraded model for $L$ (embedded in the 
canonical Jacobi resolution \ $\Uds\to G$) \ by modifying the free resolution \
$\Cds\to L$ \ of Step II, as follows\vsm :

For each \ $n\geq 0$ \ define the Jacobi algebras \ $\Jus{0}{n}$ \ 
($n=0,1,\dots$) \ to be the coproducts \ 
$\Jus{0}{0}=\Spp{1}{\lra{a}}\amalg\Spp{1}{\lra{b}}\amalg\Spp{2}{\lra{c}}$, \ 
$\Jus{0}{1}=\Spp{2}{x^{(0)}}\amalg \Spp{6}{w^{(0)}}$, \ 
$\Jus{0}{2}=\Spp{5}{y^{(0)}}$, \ $\Jus{0}{3}=\Spp{6}{z^{(0)}}$, \ and so on.

The face maps \ $d_{i}$ \ ($i=0,1,\dotsc,n$) \ are defined 
as follows: \ write \ $\lra{x}\in F(B)$ \ for the generator corresponding to an
element \ $x\in B_{\ast}'$, \ then recursively a typical DGL generator for \
$W_{n}=W_{n,\ast}$ \ (in the canonical Stover resolution \ $\Wds(B)$,
defined in \cite[\S 5.3]{BlaHR}) \ is \
$\lra{\alpha}$, \ for \ $\alpha\in W_{n-1}$, \ so an
element of \ $W_{n}$ \ is a sum of iterated Lie products of elements 
of \ $B'_{\ast}$, \ arranged within \ $n+1$ \ nested pairs of brackets \
$\lra{\lra{\dotsb}}$. \ With this notation, \ the $i$-th face map of \
$\Wds$ \ is ``omit $i$-th pair of brackets'', and the
$j$-th degeneracy map is ``repeat $j$-th pair of brackets''. \ We
assume the bracket operation \ $\lra{-}$ \ is linear \ -- \ i.e., that \ 
$\lra{\alpha x+\beta y}=\alpha\lra{x}+\beta\lra{y}$ \ for \
$\alpha,\beta\in\QQ$ \ and \ $x,y\in B$\vsm .

%
\noindent \textbf{1}. \ For \ 
$y^{(0)}=\lra{3[\lra{[\lra{a},\lra{a}]},\lra{\lra{a}}]}\in \Ju{0}{2}{3}$, \ 
we have \ $d_{0}(y^{(0)})=3[\lra{[\lra{a},\lra{a}]},\lra{\lra{a}}]=
3[x^{(0)},a^{(0)}]$ \ and \
$d_{2}(y^{(0)})=\lra{3[\lra{[a,a]},\lra{a}]}$, \ while \ $d_{1}(y^{(0)})=
\lra{3[[\lra{a},\lra{a}],\lra{a}]}$. \  This no longer vanishes as in
Step II, since the Jacobi identity does not hold in \ $\DGN$, \ but we have
an element \ $y^{(1;1)}\DEF
\lra{\lamt{\lra{a}}{\lra{a}}{\lra{a}}}\in \Ju{1}{1}{4}$ \ 
(in the notation of \cite[Example 7.8]{BlaHR}, \ with \ 
$\partial(y^{(1;1)})=d_{1}(y^{(0)})$.

On the other hand, we also have an element \ 
$y^{(1;2)}\DEF\lra{3[\lra{x},\lra{a}]}\in \Ju{1}{1}{4}$ \ (which we
denoted simply by \ $y^{(1)}$ \ in Step II, \ with \ 
$\partial(y^{(1;2)})=d_{2}(y^{(0)})$. \ The simplicial identity \ 
$d_{1}d_{2}=d_{1}d_{1}$ \ implies that \ $d_{1}(y^{(1;2)})-d_{1}(y^{(1;1)})=
\lra{3[x,a]-\lamt{a}{a}{a}}$ \ is a \
$\partial_{J}$-cycle, so we have \ $y^{(2;1,2)}\in \Ju{2}{0}{5}$ \ with \ 
$\partial_{J}(y^{(2;1,2)})=\lra{3[x,a]-\lamt{a}{a}{a}}$\vsm.

%
\noindent \textbf{2}. \ For \ $z^{(0)}=
\lra{4[\lra{3[\lra{[\lra{a},\lra{a}]}, \lra{\lra{a}}]},\lra{\lra{\lra{a}}}] + 
6[\lra{[\lra{\lra{a}},\lra{\lra{a}}]},\lra{\lra{[\lra{a},\lra{a}]}}]}$ \ 
in \ $\Ju{0}{3}{4}$\vsm :

\begin{enumerate}
\renewcommand{\labelenumi}{(\alph{enumi})~}
\item $z^{(1;1)}\DEF
\lra{\lambda_{3}(x^{(0)}\otimes \lra{\lra{a}}\otimes\lra{\lra{a}})}$, \ with
$$
d_{1}(z^{(0)}) = \lra{6(2[[x^{(0)},\lra{\lra{a}}],\lra{\lra{a}}],
+ [[\lra{\lra{a}},\lra{\lra{a}}],x^{(0)}])}=\partial_{J}(z^{(1;1)}),
$$
\item 
$z^{(1;2)}\DEF\lra{4[\lamt{\lra{a}}{\lra{a}}{\lra{a}},\lra{\lra{a}}]}$, \ 
with \ 
$$
d_{2}(z^{(0)}) = 
\lra{4[[[a^{(0)},a^{(0)}],a^{(0)}],\lra{\lra{a}}]}=\partial_{J}(z^{(1;2)})
\text{ \ \ since \ \ }[x^{(0)},x^{(0)}]=0).
$$
\item $z^{(1;3)}=\lra{4[\lra{3[\lra{[a,a]}, \lra{a}]},\lra{\lra{a}}] + 
6[\lra{[\lra{a},\lra{a}]},\lra{\lra{[a,a]}}]}\in \Ju{1}{2}{5}$, \ with \ 
$$
d_{3}(z^{(0)})=\partial_{J}(z^{(1;3)})\vsm .
$$
\end{enumerate}

%
\noindent \textbf{3}. \ Next, in \ $\Juds{2}$\vsm :
\begin{enumerate}
\renewcommand{\labelenumi}{(\alph{enumi})}
%
\item The simplicial identity \ $d_{1}d_{1}=d_{1}d_{2}$ \ implies that \ 
$$
d_{1}(z^{(1;1)})-d_{1}(z^{(1;2)})=
\lra{4[\lamt{\lra{a}}{\lra{a}}{\lra{a}},\lra{a}]-
6\lamt{[\lra{a},\lra{a}]}{\lra{a}}{\lra{a}}}
$$
is a \ $\partial_{J}$-cycle \ -- \ and indeed we have \ 
$z^{(2;1,2)}\DEF
\lra{\lamf{\lra{a}}{\lra{a}}{\lra{a}}{\lra{a}}}\in \Ju{2}{2}{6}$ \ with \ 
$\partial(z^{(2;1,2)})=d_{1}(z^{(1;1)})-d_{1}(z^{(1;2)})$.
%
\item \ Similarly, 
$$
d_{1}(z^{(1;3)})-d_{2}(z^{(1;1)})=
\lra{12[[\lra{x},\lra{a}],\lra{a}] + 6[[\lra{a},\lra{a}],\lra{x}]-
6\lamt{\lra{[a,a]}}{\lra{a}}{\lra{a}}}
$$ 
is a \ $\partial_{J}$-cycle, so we have \ $z^{(2;1,3)}\DEF
\lra{6\lamt{\lra{x}}{\lra{a}}{\lra{a}}}\in \Ju{2}{2}{6}$ \ 
with \ 
$\partial(z^{(2;1,3)})=d_{1}(z^{(1;3)})-d_{2}(z^{(1;1)})$.
\item 
$d_{2}(z^{(1;3)})-d_{2}(z^{(1;2)})=
\lra{4[\lra{3[x,a]},\lra{a}]+6[\lra{[a,a]},\lra{x}]-
4[\lra{\lamt{a}{a}{a}},\lra{a}]}$ \ is a $\partial_{J}$-cycle, hit by \
$z^{(2;2,3)}\DEF\lra{4[\lra{y},\lra{a}]+3[\lra{x},\lra{x}]}\in
\Ju{2}{2}{6}$\vsm .
\end{enumerate}

%
\noindent \textbf{4}. \ Finally, we have \ 

\begin{enumerate}
\renewcommand{\labelenumi}{(\alph{enumi})}
\item \ For \ $d_{1}z^{(2;1,2)}=\lra{\lamf{a}{a}{a}{a}}$ \ we have \ 
$$
\partial_{J}(d_{1}z^{(2;1,2)})=
d_{1}d_{2}(z^{(1;2)})-d_{1}d_{2}(z^{(1;1)})=
\lra{4[\lamt{a}{a}{a},a]-6\lamt{[a,a]}{a}{a}}
$$
\item  For \ $d_{1}z^{(2;1,3)}=\lra{6\lamt{x}{a}{a}}$ \ we have \ 
$$
\partial_{J}(d_{1}z^{(2;1,3)})=
d_{1}d_{2}(z^{(1;1)})-d_{1}d_{1}(z^{(1;3)})=
\lra{6\lamt{[a,a]}{a}{a}-12[[x,a],a]-6[[a,a],x]}
$$
\item For \ $d_{1}z^{(2;2,3)}=\lra{4[y,a]+3[x,a]}$ \ we have \ 
$$
\partial_{J}(d_{1}z^{(2;2,3)})=
d_{1}d_{2}(z^{(1;2)})-d_{1}d_{1}(z^{(1;3)})=
\lra{4[\lamt{a}{a}{a},a]-12[[x,a],a]-6[[a,a],x]},
$$
\end{enumerate}

So there is an element  \ $z^{(3;1,2,3)}\in \Ju{3}{0}{7}$ \ with 
$$
\partial_{J}(z^{(3;1,2,3)})=
\lra{\lamf{a}{a}{a}{a} - 6\lamt{x}{a}{a} + 4[y,a]+3[x,a]}.
$$
This defines the second-order homotopy operation we are interested in
(which also has an integral version).
%
%
\sect{Refined moduli spaces}
\label{crms}

Since neither the cohomolgy algebra nor the homotopy Lie algebra of a
rational space determine the other, we can refine the coarse partition
of the rational moduli space $\M$ by specifying both: 

\begin{defn}\label{drefmod}\stepcounter{subsection}
For each $1$-connected graded-commutative algebra $H$ and graded Lie
algebra $P$ over $\QQ$, \ let \ $\M_{H,P}$ \
denote the collection of all simply-connected rational homotopy types
of spaces $X$ with \ $H^{\ast}(X;\QQ)\cong H$ \ and \
$\pis(X)\otimes\QQ\cong P$. 
\end{defn}

Note that we cannot simply identify this as a quotient variety with \
$\M_{H}\cap\M_{P}$ \ in any natural way, since points 
in \ $\M_{H}$ \ are represented by filtered DGA models, while those in \ 
$\M_{P}$ \ are represented by filtered DGA models.

The set \ $\M_{H,P}$ \ may, of course, be empty; but Lemaire and Sigrist have
shown that it can also be infinite (see \cite[\S 3]{LSigD}). In order
to analyze this refined moduli space, we shall need one additional
ingredient.

\subsection{The Quillen equivalences}
\label{sqe}\stepcounter{thm}

The fact that \ $\ho\DGLz$ \ and \ $\ho\DGAo$ \ are equivalent to
each other, leads us to expect a direct algebraic relationship between
the corresponding model categories, which in fact exists, and 
may be described as follows\vsm :

For simplicity, restrict attention to simply-connected spaces of
finite type. We may therefore assume that \ $A^{i}_{X}$ \ is finite dimensional
for each \ $i\geq 0$ \ (and of course \ $A^{0}_{X}=\QQ$, \ $A^{1}_{X}=0$). \ 
Taking the vector space dual of \ $A^{\ast}_{X}$, \ we obtain a
$1$-connected differential graded-cocommutative coalgebra \
$(C_{\ast}^{X},\delta)$, \ whose homology is \ $H_{\ast}(X;\QQ)$, \
with the usual coalgebra structure (dual to the ring structure in
cohomology). Let \ $\DGC$ \ denote the category of such coalgebras.

Quillen defined a pair of adjoint functors 
$$
\DGLz~\adj{\C}{\LL}~\DGC~
$$
as follows\vs:

\noindent\textbf{I:} \ \ \ Given a DGC \ $(C_{\ast},d)$, \ let \ 
$\LL(C_{\ast}):=(\Prim(\Omega C_{\ast}),\partial)$ \ denote 
the graded Lie algebra of primitives in the cobar construction of \ 
$C_{\ast}$, \ constructed as follows:
 
If \ $C_{\ast}\cong \QQ\oplus \bar{C}_{\ast}$ \ (where \ 
$\bar{C}_{\ast}=C_{\geq 2}$, \ in our case), \ and \ 
$\Sigma^{-1}V$ \ is the graded vector space $V$ \ shifted downwards (so that \
$(\Sigma^{-1}V)_{i}:=V_{i+1}$, \ with \  $\sigma^{-1}v\leftrightarrow v$), \ 
let \ $\Omega C_{\ast}:=T(\Sigma^{-1}\bar{C}_{\ast})$ \ denote
the tensor algebra on \ $\Sigma^{-1}\bar{C}_{\ast}$ \ with 
$$
\partial(\sigma^{-1}x):=-\sigma^{-1}(dx)+
\frac{1}{2}\sum_{i}~(-1)^{|c'_{i}|}[\sigma^{-1}c_{i}',\sigma^{-1}c_{i}'']~,
$$
where \ $\bar{\Delta} c:=\sum_{i} c_{i}'\otimes c_{i}''$ \ is the (reduced)
comultiplication in \ $\bar{C}_{\ast}$, \ and \ $[~,~]$ \ is the commutator
in \ $T(\Sigma^{-1}\bar{C}_{\ast})$\vsm .

\noindent\textbf{II:} \ \ \ Given a DGL \ $(L_{\ast},\partial)$, \ let \ 
$\C(L_{\ast},d)$ \ be the DGC \ $(\Lambda(\Sigma L_{\ast}),d)$, \ where \ 
$\Sigma L_{\ast}$ \ is \ $L_{\ast}$ \ shifted upwards, \ $\Lambda V$ \
denotes the cofree graded coalgebra cogenerated by the
graded vector space $V$, and the coderivation $\partial$ \ is defined
by \ $d=d'-d''$, \ where
%
\begin{equation}\label{eone}
\begin{split}
d'(\sigma x_{1}\wedge\dotsc \wedge\sigma x_{n})~:=~&
\sum_{i=1}^{n}~(-1)^{i+\sum_{j<i}|x_{j}|}
\sigma x_{1}\wedge\dotsc \wedge\sigma \partial x_{i}\wedge\dotsc
\wedge\sigma x_{n}),\\
d''(\sigma x_{1}\wedge\dotsc \wedge\sigma x_{n})~:=~&
\sum_{1\leq i<j\leq n}~(-1)^{|x_{i}|}(\pm)
\sigma[x_{i},x_{j}]\wedge \sigma x_{1}\dotsc 
\widehat{\sigma x_{i}}\dotsc\widehat{\sigma x_{j}}\dotsc\sigma x_{n})~,
\end{split}
\end{equation}
and the sign \ $(\pm)$ \ is determined by \ 
$$
\sigma x_{1}\wedge\dotsc \wedge\sigma x_{n}=
(\pm)\sigma x_{i}\wedge\sigma x_{j}\wedge \sigma x_{1}\wedge\dotsc 
\widehat{\sigma x_{i}}\dotsc\widehat{\sigma x_{j}}\dotsc\wedge\sigma x_{n})~.
$$

See \cite[App.~B]{QuiR} or \cite[\S 22]{FHThR} for the details. 

\subsection{Joint deformations}
\label{sjd}\stepcounter{thm}

This suggests the following approach to the refined moduli
space for a graded Lie algebra $P$ and graded cocommutative coalgebra $H$, 
where we assume $H$ and $P$ are both of finite type, $H$ is
$n$-connected and $P$ is \ $(n-1)$-connected \ ($n\geq 1$), \ and \ 
$H_{n+1}\cong P_{n}$:

Let \ $(B,\partial)$ \ be the bigraded model for \ $(P,0)$, \ and \ 
$(C,d)=(C_{\ast}B,d'-d'')$ \ the corresponding coalgebra (\S \ref{sqe}~II). 
Note that \ $(C,d)$ \ is cofibrant, but not necessarily minimal. As we
deform \ $(B,\partial)$ \ to produce all of \ $\M_{P}$, \ we 
change only \ $d'$ \ in the coalgebra model, obtaining \
$(\hat{B},\hat{\partial})$, \ say. At the same time, we have
the deformations \ $(\tilde{A},\tilde{d})$ \ say, of the formal DGC \ $(H,0)$; \ 
and for each pair \ $\lra{\hat{B},\tilde{A}}$, \ we have the variety \ 
$V_{\lra{\hat{B},\tilde{A}}}$ \ of all DGC maps \ 
$\phi:(\hat{B},\hat{\partial})\to(\tilde{A},\tilde{d})$ \ between
them. Since both are cofree, these are determined by maps of graded
vector spaces. The condition that $\phi$ be a quasi-isomorphism
translates into a series of rank conditions, so that we get a
semi-algebraic set parametrizing all such pairs which are equipped
with a quasi-isomorphism between them \ -- \ and in particular, the
requisite homology $H$ and homotopy groups $P$.

%
%
\sect{Nerves and moduli spaces}
\label{cnms}

Since ordinary (integral) homotopy types do not have any known
differential graded models, there is no hope of generalizing the
deformation approach of sections \ref{crhdt}-\ref{crms} to cover them, too.
However, in \cite{DKanCD}, Dwyer and Kan suggested an approach to
such ``moduli problems'' based on the concept of nerves, which has
proved useful conceptually in a number of contexts.

\begin{defn}\label{dnerve}\stepcounter{subsection}
The \emph{nerve} \ $\Ner C$ \ of a small category $C$ is a simplicial
set whose $k$-simplices are the sequences of $k$ composable arrows in
$C$, with \ $d_{i}$=``delete $i$-th object and compose'', \ 
$s_{j}$=``insert identity arrow after $j$-th object''. The geometric
realization of the nerve is called the \emph{classifying space} of
$\C$, written  \ $BC:=|\Ner C|$.
\end{defn}

The nerve was originally defined by Segal in \cite{SegC}, based on ideas of
Grothendieck; Quillen, in \cite{QuiK} helped clarify the close connection
between nerves and homotopy theory, as evinced in the following properties:

\begin{enumerate}
\renewcommand{\labelenumi}{\arabic{enumi})~}
\item The functor \ $\Ner:\Cat\to\Ss$ \ takes natural transformations
to homotopies (cf.\ \cite[\S 2, Prop.\ 2]{QuiK}), so that
\item if a functor $F$ has a left or right adjoint, then \ $\Ner F$ \ 
is a homotopy equivalence; more general conditions are provided by
\cite[\S 2, Thm.\ A]{QuiK}). If $C$ has an initial or final object,
then \ $\Ner C\simeq\ast$. 
\item The nerve of a functor is not often a fibration; conditions
when \ $\Ner F$ \ fits into a quasi-fibration sequence are provided by 
\cite[\S 2, Thm.\ B]{QuiK}.
\item Similarly, \ $\Ner C$ \ is not usually a Kan complex, unless $C$
is a groupoid \ (cf.\ \cite[I, Lemma 3.5]{GJarS}). 
\end{enumerate}

\begin{defn}\label{dcc}\stepcounter{subsection}
A \emph{classification complex} in a model category $\C$ (see 
\cite[\S 2.1]{DKanCD}) is the nerve of some subcategory $\D$ of $\C$,
all of whose maps are weak 
equivalences, and which includes all weak equivalences whose source or
target is in $\D$. The category $\D$ is only required to be 
\emph{homotopically} small (cf.\ \cite[\S 2.2]{DKanF}) \ -- \ that is, \ 
$\Ner\D$ \ has a set of components, and its homotopy groups (at each
vertex) are small.  
\end{defn}

This construction has two main properies:

\begin{enumerate}
\renewcommand{\labelenumi}{\alph{enumi})~}
\item The components of \ $\Ner\D$ \ are in one-to-one
correspondence with the weak homotopy types of $\D$ in $\C$.
\item The component of \ $\Ner\D$ \ corresponding to an object \ $X\in\C$ \
is weakly homotopy equivalent to the classifying space of the monoid \
$\haut X$ \ of self weak equivalences of $X$ (cf.\ \cite[Prop.\ 2.3]{DKanCD}).
\end{enumerate}

\begin{remark}\label{rcc}\stepcounter{subsection}
Taking the model category \ $\C:=\To$ \ of connected pointed spaces,
with \ $\D=\W_{\C}$ \ the subcategory of \emph{all} weak homotopy
equivalences in $\C$, \ $\M^{\C}:=\Ner\W_{\C}$ \ is the
candidate suggested by the Dwyer-Kan approach for the moduli space of
pointed homotopy types. There is also an unpointed version,
of course; on the other hand, allowing for non-connected spaces merely
complicates the combinatorics of $\M$, without adding any new information.

One might ask in what sense this qualifies as a moduli space (aside
from having the right components). Even though the analogy with the
classical examples mentioned in the introduction is not clear-cut,
note that \ $\M^{\To}$ \ as defined using the nerve (which is natural
choice for a space to associate to a category) is related
to \ $\M^{\Ton^{\QQ}}$ \ of \S \ref{smrht}, \ and the latter does
exhibit many of the attributes listed in \S \ref{sms}. Moreover, as we
shall see in the next section, this construction can be used to
interpret the obstruction theory of \S \ref{sot}, as well as its
integral analogue, and also to describe $\M$ as the limit of a tower
of fibrations, which give increasingly accurate approximations to $\M$.  

Finally, the monoid \ $\haut X$ \ of self equivalences corresponds to the
mapping class group of self diffeomorphisms of a surface
(Sec.~\ref{cint}), whose classifying space is closely related (or even
homotopy equivalent) to the moduli space.
\end{remark}

\subsection{The model categories}
\label{stmc}\stepcounter{thm}

In order to provide a uniform treatment in different model
categories, and in particular to allow for the comparisons mentioned
in \S \ref{smr}, it is useful to consider a resolution model
category structure (\S \ref{srmc}) on a category \ $\C=s\E$ \ of
simplicial objects over another category $\E$. In fact, one can do 
this at two different levels, and in some sense the comparison between
these is the heart of this approach\vs: 

\noindent\textbf{I.} \ \ \ The \emph{topological} level \ -- \ where
$\E$ can take several forms\vsm: 

\begin{enumerate}
\renewcommand{\labelenumi}{(\alph{enumi})~}
\item The category \ $\To$ \ of connected pointed spaces.
\item The subcategory \ $\To^{\QQ}$ \ of pointed connected spaces
having rational universal covers (but arbitrary fundamental group).

Note that one can approximate any pointed connected space \ 
$X\in\To$ \ by its fibrewise rationalization \ $X'_{\QQ}$ \ (cf.\
\cite[I, 8.2]{BKaH}), \ which lies in \ $\To^{\QQ}$: \ if \
$\tilde{X}\to X\to B(\pi_{1}X)$ \ is the universal covering space 
fibration for $X$, then  \ $X'_{\QQ}$ \ fits into a functorial
fibration sequence \ $\widetilde{X_{\QQ}}\to X'_{\QQ}\to B(\pi_{1}X)$, \ 
in which \ $\widetilde{X_{\QQ}}$ \ is the usual rationalization of the universal
cover.  However, algebraic DGL, DGA, or DGC models do not generally
extend to this case (unless \ $\pi_{1}X$ \ is finite and acts
nilpotently on the higher groups \ -- \ see, e.g., \cite{TriE}). 
\item The subcategory \ $\Ton^{\QQ}$ \ of $1$-connected 
pointed rational spaces (in \ $\To^{\QQ}$). \ This can be replaced
by the Quillen equivalent model categories \ $\DGLz$ \ or \ $\DGAo$.
\item Other variants are possible \ -- \ for instance, we could consider
functor categories over $\E$ \ -- \ that is, diagrams \ $D\to\E$ \ for
a fixed small category $D$ (see \cite{BJTurnM})\vsm.
\end{enumerate}

\noindent\textbf{II.} \ \ \ The \emph{algebraic} level \ -- \ where
$\E$ is correspondingly\vsm : 

\begin{enumerate}
\renewcommand{\labelenumi}{(\alph{enumi})~}
\item The category \ $\PAlg$ \ of \emph{\Pa s} \ -- \ that is,
positively graded groups \ $G_{\ast}$, \ abelian in dimensions $\geq 2$, \ 
equipped with an action of the primary homotopy operations
(Whitehead products, compositions, and \ $G_{1}$-action) satisfying
the usual identities (see \cite{BlaA} for a more explicit definition).
\item The subcategory \ $\PAlg^{\QQ}$ \ of \emph{rational \Pa s}, which
are (positively) graded Lie algebras equipped with a ``fundamental
group action'' of an arbitrary group $\pi$, satisfying the usual
identities (see \cite{HilJ}). 
\item When \ $\pi=0$, we obtain the subcategory \ 
$\PAlg^{\QQ}_{1}\cong\gr\Lie$ \ of simply-connected rational \Pa s,
which are just graded Lie algebras over $\QQ$. 
\item There is also a concept of \Pa s for arbitrary diagrams (again,
see \cite{BJTurnM}).
\end{enumerate}

%
%
\sect{Approximating classification complexes}
\label{cacc}

The classification complex of a model category $\C$, as defined in 
\S \ref{dcc}, may appear to be a somewhat artificial marriage of a
traditional moduli space, whose set of components correspond to the
(weak) homotopy types of $\C$, and the individual components, which are
classifying spaces of the form \ $B\haut X$. \ It is not clear at
first glance why such complexes might be useful. To understand this,
we show how \ $\M^{\C}$ \ can be approximated by a tower of
fibrations, in a way that elucidates the obstruction theory of \S \ref{sot}:

\subsection{Postnikov systems and Eilenberg-Mac Lane objects}
\label{sps}\stepcounter{thm}

Recall from \cite[\S 1.2]{DKanO} that functorial Postnikov towers \ 
$$
\Xd\to\dotsc\to\tP{n}\Xd\to\tP{n-1}\Xd\to\dotsc\to\tP{0}\Xd
$$
may be defined in categories of the form \ $\C=s\E$ \ using the matching space
construction of \cite[X,\S 4.5]{BKaH}. \ By considering the fibers of
successive Postnikov sections, we can define the analogue of homotopy
groups in each such category, and find that they are corepresented, as
one might expect, by appropriate suspensions of the ``good objects''
in $\E$. For both \ $\E=\To$ \ and \ $\E=\PAlg$, \ we find that the
natural ``homotopy groups'' \ $\pin{n}{\Xd}$ \ of any \ $\Xd\in s\E$ \ (the
\emph{bigraded groups} of \cite{DKStB}) \ take values in \ $\PAlg$, \
with the obvious modifications for the rational variants. 

Note that for a simplicial space \ $\Xd\in s\To$, \ applying \ $\pis$ \ 
in each simplicial dimension yields a simplicial \Pa\ \ $\pis\Xd$, \
and the two sequences of \Pa s \ $(\pin{n}{\Xd})_{n=0}^{\infty}$ \ and \  
$(\pi_{n}\pis\Xd)_{n=0}^{\infty}$ \ fit into a ``spiral long
exact sequence'': 
%
\setcounter{equation}{\value{thm}}\stepcounter{subsection}
\begin{equation}\label{etwo}
\dotsc \pi_{n+1}\pis\Xd \xra{\partial^{\star}_{n+1}} \Omega\pin{n-1}{\Xd}
\xra{s_{n}} \pin{n}{\Xd} \xra{h_{n}} \pi_{n}\pis\Xd \xra{\partial^{\star}_{n}} 
\dotsb \pin{0}{\Xd}\xra{h_{0}}\pi_{0}\pis\Xd \to 0
\end{equation}
\setcounter{thm}{\value{equation}}
\noindent (see \cite[8.1]{DKStB}), in which each term is not only a
\Pa, but a module over \ $\pin{0}{\Xd}\cong\pi_{0}\pis\Xd$ \ under a
``fundamental group action'', for \ $n\geq 1$. \ Here \ $\Omega\Lambda$ \
denotes the abelian \Pa\ obtained from a \Pa\ $\Lambda$ \ by
re-indexing and suspending the operations (so that in particular \ 
$\Omega\pis X\cong\pis\Omega X$ \ for any space $X$).

Furthermore, one can construct \emph{classifying objects} \ 
$\bB\Lambda\in s\To$ \ for any \Pa\ $\Lambda$ \ (with \ 
$\pin{0}{\bB\Lambda}=\Lambda$ \ and \ $\pin{i}{\bB\Lambda}=0$ \ otherwise); 
\emph{Eilenberg-Mac Lane} objects \
$\bB(\Lambda,n)$; \ and \emph{twisted} Eilenberg-Mac Lane objects \ 
$\bB_{\Lambda}(M,n)$ \ for any \Pa\ $\Lambda$, \ $\Lambda$-module $M$,
and \ $n\geq 1$, \ with
$$
\pin{i}\bB_{\Lambda}(M,n)\cong\begin{cases}
\Lambda &~~\text{if \ }  i=0\\
M &~~\text{(as a $\Lambda$-module) \ if} \ i=n\\
0 &~~\text{otherwise.}
\end{cases}
$$

This is true more generally for resolution model categories,
under reasonable assumptions. To simplify the notation we denote by \ 
$K_{\Lambda}$ \ the analogous simplicial \Pa \
(with \ $\pi_{n}K_{\Lambda}\cong\Lambda$ \ for \ $n=0$, \ and $0$
otherwise), and similarly the Eilenberg-Mac Lane objects \
$K(\Lambda,n)$, \ and twisted Eilenberg-Mac Lane objects \ 
$K_{\Lambda}(M,n)$ \ in \ $s\PAlg$.
We shall use boldface in general to indicate
constructions in \ $s\To$, \ as opposed to \ $s\PAlg$. 

Finally, one can define natural $k$-invariants for the Postinkov
system of a simplicial object \ $\Xd\in s\E$, \ as in 
\cite[Prop.\ 6.4]{BDGoeR}, and these take values in appropriate
Andr\'{e}-Quillen cohomolgy groups (represented by the twisted
Eilenberg-Mac Lane objects). These groups are denoted respectively by \ 
$H^{n}(\Xd/\Lambda;M):=[\Xd,\bB_{\Lambda}(M,n)]_{\bB\Lambda}$, \ for a
simplicial space \ $\Xd$ \ equipped with a map to \ $\bB\Lambda$, \ 
and \ $H^{n}(\Gd/\Lambda;M):=[\Gd,K_{\Lambda}(M,n)]_{K_\Lambda}$, \ 
for a simplicial \Pa\ \ $\Gd$ \ equipped with a map to \ $K_{\Lambda}$.

It turns out that for any simplicial space \ $\Xd$ \ there
is a natural isomorphism 
%
\setcounter{equation}{\value{thm}}\stepcounter{subsection}
\begin{equation}\label{eseven}
H^{n}(\Xd/\Lambda;M)\xra{\cong}H^{n}(\pis\Xd/\Lambda;M)
\end{equation}
\setcounter{thm}{\value{equation}}
\noindent for every \ $n\geq 1$ \ (cf.\  \cite[Prop.\ 8.7]{BDGoeR}).

\subsection{Relating classification complexes}
\label{srcc}\stepcounter{thm}

We shall now show how when \ $\E=\To$ \ -- \ and more generally \ -- \
the classification complex \  
$\M^{\To}$ \ can be exhibited as the homotopy limit of a tower of
fibrations, where the successive fibers have a cohomological
description showing the relationship between \ 
$\pi_{0}\M^{\To}$ \ and the higher homotopy groups.
The tower in question is constructed essentially by taking successive Postnikov
sections. The idea is an old one, and is useful even in analyzing the
self-equivalences of a single space $X$ (see, for example, \cite{WilkC}). 

\begin{defn}\label{dmmodto}\stepcounter{subsection}
For a given \Pa \ $\Lambda$, we denote by \
$\D(\Lambda)=\D^{s\To}(\Lambda)$ \ the category of simplicial spaces \ 
$\Xd\in s\To$ \ such that \ $\pis\Xd\simeq\bB\Lambda$ \ (in \ $s\PAlg$) \ 
(that is, \ $\pi_{n}\pis\Xd\cong\Lambda$ \ for \ $n=0$, \ and  \ 
$\pi_{n}\pis\Xd-0$ \ otherwise). \ The nerve of \
$\D^{s\To}(\Lambda)$ \ will be denoted by \ $\M_{\Lambda}$.

The ``pointed'' version is the nerve of the category \ $\R(\Lambda)$ \ 
of pairs \ $(\Xd,\rho)$, \ where \ $\Xd\in s\To$ \ and \ 
$\rho:\bB\Lambda\to\pis\Xd$ \ is a specified weak equivalence in \
$s\PAlg$ \ (again with weak equivalences as morphisms).   
\end{defn}

Although \ $\M_{\Lambda}$ \ is the more natural object of interest in 
our context, we actually study \ $\R(\Lambda)$.  \ 
As noted in \cite[\S 1.1]{BDGoeR}, \ there is a fibration sequence \  
%
\setcounter{equation}{\value{thm}}\stepcounter{subsection}
\begin{equation}\label{ethree}
\Ner\R(\Lambda)\to\M_{\Lambda}\to B\Aut(\Lambda)~,
\end{equation}
\setcounter{thm}{\value{equation}}
\noindent where \ $\Aut(\Lambda)$ \ is the group of automorphisms of
the \Pa $\Lambda$. 

\begin{defn}\label{dnp}\stepcounter{subsection}
For each \ $n\geq 1$, \ let \ $\R_{n}(\Lambda)$ \ denote the category of
$n$-Postnikov sections under \ $K_{\Lambda}$ \ -- \ that is, the 
objects of \ $\R_{n}(\Lambda)$ \ are pairs \ $(\Xp{n},\rho)$, \ where \ 
$\Xp{n}\in s\To$ \ is a simplicial space such that \ 
$\bP{n}\Xp{n}\simeq\Xp{n}$ \ (Postnikov sections in \ $s\To$), \ and \ 
$\rho:\tP{n}K_{\Lambda}\to\tP{n}\pis\Xp{n}$ \ is a weak equivalence. 
The morphisms of \ $\R_{n}(\Lambda)$ \ are weak equivalences of simplicial spaces
compatible with the maps $\rho$ up to weak equivalence of simplicial \Pa s.
\end{defn}

\subsection{A tower of realization spaces}
\label{strs}\stepcounter{thm}

Given $\Lambda$, the Postnikov section functors \ 
$\bP{n}:s\To\to s\To$ \ of \S \ref{sps} induce compatible functors \ 
$\Phi_{n}:\R(\Lambda)\to\R_{n}(\Lambda)$ \ and \ 
$F_{n}:\R_{n+1}(\Lambda)\to\R_{n}(\Lambda)$; \ and as in 
\cite[Thm.\ 9.4]{BDGoeR} and \cite[Thm.\ 3.4]{DKanCD}, these in turn 
induce a weak equivalence 
%
\setcounter{equation}{\value{thm}}\stepcounter{subsection}
\begin{equation}\label{efour}
\Ner\R(\Lambda)\to \holim_{n} \Ner\R_{n}(\Lambda).
\end{equation}
\setcounter{thm}{\value{equation}}
\noindent Combining \eqref{ethree} and \eqref{efour}, we may try to obtain 
information about the space of realizations \ $\M_{\Lambda}$ \ by studying the 
successive stages in the tower
%
\setcounter{equation}{\value{thm}}\stepcounter{subsection}
\begin{equation}\label{efive}
\dotsc \Ner\R_{n+1}(\Lambda) \xra{\Ner F_{n}}\Ner\R_{n}(\Lambda) 
\xra{\Ner F_{n-1}}\dotsc\to\Ner\R_{1}(\Lambda).
\end{equation}
\setcounter{thm}{\value{equation}}

%
%
\sect{Analyzing the tower}
\label{catt}

The first step in analyzing the tower \eqref{efive} is to
understand when the successive fibers are non-empty, and if
so, to count their components. 
In the rational case, too, our main task was identifying 
the components of the space of rational homotopy types, and a partial
ordering on the components was induced by the successive deformations.
The problem of empty fibers did not arise there, since all DGLs (or DGAs)
are realizable by rational spaces. 
However, the fact that we have a ordered the successive choices in a 
tower, rather than a tree as in Theorem \ref{ttwo}, suggests that we
can describe them by means of an obstruction theory, as follows\vsm :

Assume given a point  \ $(\Xp{n},\rho)$ \ in \ $\Ner\R_{n}(\Lambda)$, \ 
so that \ $\Xp{n}\in s\To$ \ is a simplicial space which is an
$n$-Postnikov stage (for some putative simplicial space $Y$
realizing the given \Pa\ $\Lambda$), \ and \
$\rho:\bK_{\Lambda}\to\tP{n}\pis\Xp{n}$ \ is a choice of a  
weak equivalence. We can then reinterpret \cite[Prop.\ 9.11]{BDGoeR}
as saying:

%
%
\begin{prop}\label{ptwo}\stepcounter{subsection}
$\Xp{n}\in\R_{n}(\Lambda)$ \ extends to an \ $(n+1)$-Postnikov stage in \ 
$\R_{n+1}(\Lambda)$ \ if an only if the \ $n+1$-st $k$-invariant for
the simplicial \Pa\ \ $\pis\Xd$ \ vanishes in \
$H^{n+3}(\Lambda;\Omega^{n+1}\Lambda)$. 
\end{prop}

The spiral long exact sequence \eqref{etwo} actually determines the
homotopy ``groups'' of the simplicial \Pa\ \ $\pis\Xp{n}$ \ completely:
%
\setcounter{equation}{\value{thm}}\stepcounter{subsection}
\begin{equation}\label{esix}
\pi_{k}\pis\Xp{n}\cong\begin{cases}
        \Lambda & \text{for \ } k=0,\\
        \Omega^{n+1}\Lambda & \text{for \ }k=n+2\\
         0 & \text{otherwise}.
\end{cases}
\end{equation}
\setcounter{thm}{\value{equation}}

However, \eqref{esix} in itself does not imply that \ $\pis\Xp{n}$ \
is an twisted Eilenberg-Mac Lane object \
$K_{\Lambda}(\Omega^{n+1}\Lambda,n+2)$ \ in \ $s\PAlg$ \ -- \ for that
to happen, the map \ $\pis\Xp{n}\to K_{\Lambda}$ \ must have a
section (whose existence is equivalent to the vanishing of the
$k$-invariant in Proposition \ref{ptwo}).

\begin{remark}\label{robst}\stepcounter{subsection}
It may help to understand why if one considers a simplicial group \ 
$K\in\G$ \ (as a model for a connected topological space). The
fundamental group then appears as \ $\Gamma=\pi_{0}K$, \ the set of
path components of $K$ (all homotopy equivalent to each other) \ -- \ 
which happens to have a group structure. When \ 
$K=K^{\Gamma}(M,n)=K(M,n)\ltimes\Gamma$ \ is a twisted Eilenberg-Mac Lane 
object, the choice of the section \ $s:B\Gamma\to K$ \ is what
distinguishes $K$ from a disjoint union (of cardinality \ $|\Gamma|$) \ 
of copies of ordinary Eilenberg-Mac Lane objects \ $K(M,n)$ \ -- \ that is, \ 
$L:=\coprod_{|\Gamma|}~K(M,n)$. \ Even though we can put a
group structure on \ $\pi_{0}L$ \ so as to make it abstractly
isomorphic to $\Gamma$, we will not get the right action of $\Gamma$ on
$M$ (which, in the case of $K$, appears in the usual way by
conjugation with any representative of \ $\gamma\in\Gamma=\pi_{0}K$). 
\end{remark}

\subsection{Distinguishing between liftings}
\label{sdl}\stepcounter{thm}

Since the weak homotopy types of the realizations of a given \Pa\
$\Lambda$ (that is, components of \ $\M^{\To}_{\Lambda}$)
are in one-to-one correspondence with the weak homotopy types of the 
simplicial spaces $Y$ realizing \ $K_{\Lambda}$ \ in \ $s\To$ \
(if any), we should be able to describe them inductively in terms of
the simplicial-space $k$-invariants of $Y$, and thus of the
successive Postnikov approximations \ $\Xp{n}$.  

So, one naturally expects that there will be a correspondence between
the sections \ $s_{n}$ \ and the $k$-invariants \ $\bk{n}$ \ in \ $s\To$. \ 
However, by \cite[Prop 9.11]{BDGoeR} (using the identification of
\eqref{eseven}) we know that the map \ 
$\phi:K_{\Lambda}(\Omega^{n+1}\Lambda,n+2)
\to K_{\Lambda}(\Omega^{n+1}\Lambda,n+2)$ \ 
corresponding to \ $\bk{n}$ \ must be a weak equivalence, in our case \ -- \ 
that is, it is within the indeterminacy of the $\bk{}$-invariants,
which is certainly contained in \ $\Aut_{\Lambda}(\Omega^{n+1}\Lambda)$. \ 

Thus, the relevant information for constructing successive 
stages in a Postnikov tower for the simplicial space $Y$ with \ 
$\pis Y\simeq K_{\Lambda}$ \
(if it exists), aside from the \Pa\ $\Lambda$, is \emph{not} its 
$\bk{}$-invariants, nor the $k$-invariants for the 
simplicial \Pa s \ $\pis\bP{n}\bB\Lambda$ \ ($n=1,2,\dotsc$). \
Instead, it is the seemingly innocuous choice of the section \ 
$s_{n}:K_{\Lambda}\to K_{\Lambda}(\Omega^{n+1}\Lambda,n+2)$, \ with
which any twisted Eilenberg-Mac~Lane object (of simplicial spaces or \Pa s) is
equipped:

%
%
\begin{prop}\label{pthree}\stepcounter{subsection}
If \ $\Xp{n}\in\R_{n}(\Lambda)$ \ can be lifted to \
$\R_{n}(\Lambda)$, \ the possible choices for the \ $(n+1)$-st
Postnikov  stage \ $\Xp{n+1}$ \ are determined by the choices of
sections \ $s_{n}:K_{\Lambda}\to\pis\Xp{n}$, \ which correspond to
elements in \ $H^{n+2}(\Lambda;\Omega^{n+1}\Lambda)$.
\end{prop}

\subsection{The components}
\label{scomp}\stepcounter{thm}

As noted above, the Dwyer-Kan concept of classification complexes has
one advantage over the traditional moduli spaces, in that the topology
of each component encodes further information about the corresponding
homotopy type: namely, the topological group of self-equivalences. 

Such groups are generally hard to analyze, as one might expect from the
topological analogue (see, e.g., \cite{RutSS}). 
One advantage of the tower \eqref{efive} is that the successive fibers
are generalized Eilenberg Mac-Lane spaces:

%
%
\begin{prop}[\protect{\cite[Prop.~9.6\emph{ff.}]{BDGoeR}}]
\label{pfour}\stepcounter{subsection}
For any \Pa\ $\Lambda$ and \ $n\geq 0$, \ the fiber of \ 
$\R_{n+1}(\Lambda)\to\R_{n}(\Lambda)$ \ is either empty, or has all
components weakly equivalent to \ 
$\prod_{i=1}^{n+1}~K(H^{n+3-i}_{\Lambda}(\Lambda,\Omega^{n+1}\Lambda),i)$.
\end{prop}

%
%
\sect{Comparing the rational and integral versions}
\label{ccriv}

The constructions of the rational moduli space as an
infinite-dimensional algebraic variety over $\QQ$, and of the integral
moduli space as a classification complex, are not related in any
evident way. However, given a \Pa\ \ $\Lambda\in\PAlg$ \ and 
its rationalization \ $\Lambda^{\QQ}\in\PAlg^{\QQ}$ \ (defined \ 
$\Lambda^{\QQ}_{1}:=\Lambda_{1}$ \ and \ 
$\Lambda^{\QQ}_{i}:=\Lambda_{i}\otimes\QQ$ \ for \ $i\geq 2$), \ there
is a connection between the set of components of \ $\M_{\Lambda^{\QQ}}$ \ 
and that of \ $\M_{\Lambda}$, \ and also between the respective towers
\eqref{efive}.  

\subsection{The preimage of rationalization}
\label{spir}\stepcounter{thm}

Clearly, $\Lambda$ is realizable only if \ $\Lambda^{\QQ}$ \ is
(which is automatic in the simply-connected case), but the converse
need not hold (as the example of a non-realizable torsion \Pa\ in
\cite[Prop.\ 4.3.6]{BlaO} shows). Therefore, the number of components
of \ $\M_{\Lambda^{\QQ}}$ \ need have no relation to that of \ 
$\M_{\Lambda}$, \ since there seems to be no known method for
determining all the (integral) homotopy types $X$ whose fibrewise
rationalizations \ $X'_{\QQ}$ \ (\S \ref{stmc}I(b)) are homotopy equivalent.
However, the obstruction theories of \S \ref{sot} and
Propositions \ref{ptwo}-\ref{pthree}, respectively, do provide a
framework for studying the preimage of the (fibrewise) rationalization:

Given a rationalized \Pa\ $A$, we have
\begin{enumerate}
\renewcommand{\labelenumi}{(\roman{enumi})~}
\item an algebraic question: \ which \Pa s $\Lambda$ have \ 
$\Lambda^{\QQ}\cong A$? 
\item a homotopy-theoretic question: \ what are the components of \
$\M_{\Lambda}$ \ (though different components can map to the same
component of \ $\M_{A}$ \ under the fibrewise rationalization functor)?
\end{enumerate}

For the first question, restrict attention to the simply-connected
case, so that $A$ is just a (connected) graded Lie algebra over
$\QQ$. In this case we must first consider which connected 
graded Lie algebras over $\ZZ$ have \ $L\otimes\QQ\cong A$. Given such
an $L$, we would like to classify all possible \Pa\ structures on $L$;
one possible approach to this problem is given by the obstruction theory of 
\cite[\S 12]{BPescF}, in terms of certain relative cohomology groups
(where ``relative'' involves comparing the same object in different
categories via a forgetful functor).

Theoretically, the obstruction theory of Section \ref{catt} provides
an answer to the second question, although it does not tell us how to
identify all integral homotopy types having weakly equivalent
rationalizations.

\subsection{Comparing cohomology theories}
\label{scct}\stepcounter{thm}

In view of the obstruction theories of Sections \ref{crhdt} and \ref{catt}, 
a first step towards understanding the relation between the integral
and rational classifications is to compare the cohomology theories that
house the respective obstructions. There are two main cases to consider\vs :

\noindent\textbf{I.}\hs We can compare the cohomology of a \Pa\ $\Lambda$, with
coefficients in some $\Lambda$-module $M$ (say, \
$M=\Omega^{n}\Lambda$) \ with the cohomology of the rationalization \
$\Lambda^{\QQ}$, \ with rationalized coefficients in \ $M^{\QQ}$ \
(which is a $\Lambda$-module, so in particular a \ $\Lambda^{\QQ}$-module).

In this case the short exact sequence of $\Lambda$-modules \ 
$0\to K\xra{i} M\xra{q} M^{\QQ}\to 0$ \ (for \ $K:=\Ker(q)$) \ yields
a long exact sequence in cohomology:
$$
\dotsc H^{n}(\Lambda,K)~\xra{i_{\ast}}~H^{n}(\Lambda,M)~\xra{q_{\ast}}~ 
H^{n}(\Lambda,M^{\QQ})~\xra{\delta_{\ast}}~H^{n+1}(\Lambda,K)\dotsc
$$
\noindent in the usual way. Moreover, the functors \ 
$$
\PAlg~\xra{T}~\PAlg^{\QQ}~\xra{S}\Abgp,
$$
\noindent where \ $T(-):=(-)^{\QQ}$ \ and \
$S(-):=\Hom_{\Lambda}(-,M^{\QQ})$, \ satisfy the 
conditions of \cite[Thm.~4.4]{BStoG}, since a rationalized free \Pa\
is free (and so $H$-acyclic) in \ $\PAlg^{\QQ}$. \ Thus we have a
generalized Grothendieck spectral sequence with 
%
\setcounter{equation}{\value{thm}}\stepcounter{subsection}
\begin{equation}\label{eeight}
E^{2}_{s,t}=(L_{s}\bar{S}_{t})(L_{\ast}T)\Lambda~\Rightarrow~
H^{s+t}(\Lambda^{\QQ};M^{\QQ})
\end{equation}
\setcounter{thm}{\value{equation}}
\noindent converging to the cohomology in the category \
$\PAlg^{\QQ}$). \ Here \ $L_{s}$ \ denotes the $s$-th left derived
functor, and (for simply-connected $\Lambda$) \ 
$\bar{S}:\bg\Lie\to\gr\Abgp$ \ is the functor induced by $S$, 
which exists because the homotopy groups of
any simplicial graded Lie algebra over $\QQ$ actually take value in
the category \ $\bg\Lie$ \ of bigraded Lie algebras.

For general \ $\Lambda\in\PAlg$, \ by \cite[Prop.~3.2.3]{BStoG}
we have instead a functor \ $\bar{S}:(\PAlg^{\QQ})\text{-}\PAlg\to\gr\Abgp$ \ 
whose domain is the analogue for \ $\PAlg^{\QQ}$ \ of the category of
\Pa s for spaces. Its objects are bigraded groups, endowed with an
action of all primary homotopy operations which exist for the homotopy
groups of a simplicial rational \Pa. As in the simply-connected case,
these include the bigraded Lie bracket mentioned above, and presumably
others\vs. 

\noindent\textbf{II.}\hs Alternatively, one could start with a
rationalized \Pa\ $A$ \ -- \ for simplicity, a graded Lie algebra \ -- \ 
and try to compare its cohomology (with coefficients in an
$A$-module $M$) taken in the category of Lie algebras with that
obtained by thinking of it as an ordinary \Pa. This means taking the
derived functors of the composite of 
$$
\PAlg^{\QQ}~\xra{I}~\PAlg~\xra{S}\Abgp,
$$
\noindent where $I$ is the inclusion, and \ 
$S(-):=\Hom_{\PAlg/A}(-,M)=\Hom_{\PAlg/I(A)}(-,M)$.

These do not satisfy the conditions of \cite[Thm.~4.4]{BStoG}, since a
rationalized free \Pa\ is not a free \Pa. However, we can take a
simplicial resolution \ $\Vd\to A$ \ in the category \ $s\PAlg^{\QQ}$ \ 
of simplicial graded Lie algebras, and then resolve each \ $V_{n}$ \ 
functorially in \ $s\PAlg$ \ to get a bisimplicial free \Pa\ \ $\Wdd$, \ 
with \ $Xd:=\diag\Wdd$ \ a free \Pa\ resolution of \ $I(A)$. \ If \ 
$\Ab(-)$ \ denotes the abelianzation in the over category \
$\PAlg/A$, \ then \ $\Ab\Wdd$ \ is a bisimplicial abelian object,
or equivalently, a double complex in \ $(\PAlg/A)_{\ab}$, \ and the
bicosimplicial abelian group \ $S\Wdd$ \ is \
$\Hom_{\PAlg/A}(\Ab\Wdd,M)$. \ Since 
\begin{equation*}
\begin{split}
\Hom_{\PAlg/A}(\Xd,M)~=~&\Hom_{(\PAlg/A)_{\ab}}(\diag\Ab\Wdd,M)\\
\simeq~&\Tot\Hom_{(\PAlg/A)_{\ab}}(\Ab\Wdd,M)=\Tot\Hom_{\PAlg/A}(\Wdd,M)
\end{split}
\end{equation*}
\noindent by the generalized Eilenberg-Zilber Theorem, we see that
both the spectral sequences for the bicomplex \ $S\Wdd$ \ converge to
the cohomology of \ $I(A)$. \ Since by definition \ 
$$
\pi_{v}^{t}SW_{n\bullet}\cong H^{t}(IQ_{n},M)=H^{t}(L_{n}I(-),M)(A),
$$
\noindent we obtain a cohomological spectral sequence with
%
\setcounter{equation}{\value{thm}}\stepcounter{subsection}
\begin{equation}\label{enine}
E_{2}^{s,t}=(\pi^{s}H^{t}(-,M))(L_{\ast}I)A~\Rightarrow~H^{s+t}(I(A),M).
\end{equation}
\setcounter{thm}{\value{equation}}

This is less useful than \eqref{eeight}, since we cannot identify the \ 
$E_{2}$-term explicitly as a derived functor.

\begin{remark}\label{rirat}\stepcounter{subsection}
As was pointed out in Section \ref{cmhho}, another way to relate the
integral and rational moduli spaces is geometrically, using higher
homotopy operations. This was the motivations behind \cite{BlaHR} (in
conjunction with \cite{BlaHH}). The main difficulty in using homotopy
operations in any systematic way for such a purpose is the lack of an
appropriate taxonomy. The most promising way to overcome this would be
by establishing a clear correspondence between higher operations and 
suitable cohomology groups.
\end{remark}

\end{document}